\documentclass[11pt,reqno]{amsart}

\usepackage{amsmath,amscd}
\usepackage{booktabs}
\usepackage[pdftex]{graphicx}

\usepackage[top=20mm, bottom=20mm, left=20mm, right=20mm]{geometry}
\usepackage{mathrsfs}
\usepackage{enumerate,amssymb,amsthm,multirow}

\usepackage{comment}
\usepackage{wrapfig}
\usepackage{subfig}
\usepackage[rgb]{xcolor}
\usepackage{tikz}
\usepackage{bm}

\numberwithin{equation}{section}

\newtheoremstyle{standard}
  {3pt}           
  {3pt}           
  {\itshape}      
  {\parindent}    
  {}              
  {.}             
  {.5em}          
  {}              

\providecommand{\e}[1]{\ensuremath{\times 10^{#1}}}

\theoremstyle{standard}

\newtheorem{lemma}{LEMMA}[section]

\newtheorem{remark}{Remark}[section]

\newcommand{\dblint}{\int\hskip -7pt \int}

\newcommand{\lsp}{\vspace{3mm}}

\newcommand{\vtwo}[2]{\left[\begin{array}{c} #1 \\ #2 \end{array}\right]}

\newcommand{\mtwo}[4]{\left[\begin{array}{cc}    #1 & #2 \\ #3 & #4  \end{array}\right]}

\newcommand{\vct}[1]{\bm{#1}}
\newcommand{\mtx}[1]{\mathsf{#1}}
\def\XXint#1#2#3{{\setbox0=\hbox{$#1{#2#3}{\int}$}
     \vcenter{\hbox{$#2#3$}}\kern-.5\wd0}}

\newcommand{\mA}{\mtx{A}}
\newcommand{\mtA}{\tilde{\mtx{A}}}

\newcommand{\mS}{\mtx{S}}
\newcommand{\mU}{\mtx{U}}
\newcommand{\mV}{\mtx{V}}
\newcommand{\mhU}{\hat{\mtx{U}}}
\newcommand{\mhV}{\hat{\mtx{V}}}

\newcommand{\uu}{\vct{u}}
\newcommand{\tuu}{\tilde{\vct{u}}}
\newcommand{\vv}{\vct{v}}
\newcommand{\tvv}{\tilde{\vct{v}}}

\title[Integral equations on curved surfaces]
{
A high-order accurate accelerated direct solver for acoustic scattering
from surfaces
}

\date{\today}

\begin{document}

\author{James Bremer$^{\ddagger,\S}$}
\author{Adrianna Gillman$^{*}$}
\author{Per-Gunnar Martinsson$^\dagger$}
\footnotetext[1]{Department of Mathematics, Dartmouth College}
\footnotetext[2]{Department of Applied Mathematics, University of Colorado, Boulder}
\footnotetext[3]{Department of Mathematics, University of California, Davis.}
\footnotetext[4]{Corresponding author.  \emph{E-mail address}: bremer@math.ucdavis.edu.}

\begin{abstract}
We describe an accelerated direct solver for the integral equations which
model acoustic scattering from curved surfaces.
Surfaces are specified via a collection of smooth parameterizations given on triangles,
a setting
which generalizes the typical one of triangulated surfaces,
and the integral equations are discretized via a high-order Nystr\"om method.
This allows for rapid convergence in cases in which high-order surface information is available.
The high-order discretization technique is coupled with a direct solver based on
the recursive construction of scattering matrices.  The result is a solver
which often attains $O(N^{1.5})$ complexity in the number of discretization nodes $N$ and
which is resistant to many of the pathologies which stymie iterative solvers
in the numerical simulation of scattering.
The performance of the algorithm is illustrated with numerical experiments
which involve the simulation of scattering from a variety of domains,
including one consisting of a collection of $1000$ ellipsoids
with randomly oriented semiaxes arranged in a grid, and a domain
whose boundary has $12$ curved edges and $8$ corner points.

\end{abstract}

\maketitle

\begin{section}{Introduction}


The manuscript describes techniques based on boundary integral equation
formulations for numerically solving certain
linear elliptic boundary value problems associated with acoustic scattering
in three dimensions.
There are three principal components to techniques of this type:
\vskip 1em

\begin{enumerate}[1. ]
\item
{\it Reformulation.}
The boundary value problem is first reformulated as a
system of integral equations, preferably involving operators which
are well-conditioned when considered on spaces of integrable functions.
\vskip 1em

\item
{\it Discretization.}
The resulting integral operators are then discretized.  High-order discretization
methods are to be preferred to low-order schemes and the geometry of the domain
under consideration is the principal consideration in producing high-order
convergent approaches.
\vskip 1em

\item
{\it Accelerated Solver.}
The large dense system of linear algebraic equations produced through discretization
must be solved.  Accelerated solvers which exploit the mathematical properties of the
underlying physical problems are required to solve most problems of interest.
\vskip 1em
\end{enumerate}

For boundary value problems given on two-dimensional domains
--- which result in integral operators on planar curves ---
satisfactory approaches to each of these three tasks have been described in the literature.
Uniquely solvable integral equation formulations of boundary value problems for Laplace's equation
and the Helmholtz equation are available \cite{Brakhage-Werner,Taskinen-Yla-Oijala,Barnett-Greengard,Hackbusch,Colton-Kress2};
many mechanisms for the high-order discretization of singular integral
operators on planar curves are known \cite{Colton-Kress2,Kress} and even
operators given on planar curves with corner singularities can be discretized
efficiency and to near machine precision accuracy
\cite{Bremer1,Bremer2,Bremer-Rokhlin-Sammis, Helsing-Ojala1,Helsing1};
finally, fast direct solvers which often have run times which are asymptotically
optimal in the number of discretization nodes $N$ have been constructed
\cite{Martinsson-Rokhlin, 2012_martinsson_FDS_survey,2010_gu_xia_HSS,2004_borm_hackbusch,2003_hackbusch,2003_gu,2002_hackbusch_H2,2012_ho_greengard_fastdirect}.

For three-dimensional domains, the state of affairs is much less satisfactory.
In this case, the integral operators of interest are given on curved surfaces
which greatly complicates problems (2) and (3) above.
Problem (2) requires the evaluation of large numbers of singular integrals
given on curved regions.    Standard mechanisms for addressing this problem,
like the Duffy transform followed by Gaussian quadrature,
are limited in their applicability (see \cite{Bremer-Gimbutas1} for
a discussion).  Most fast solvers (designed to address problem (3)) are based on
iterative methods and can fare poorly when the surface under consideration
exhibits complicated geometry. For example, a geometry consisting of a
$2\times 2\times 2$ grid of ellipsoids ($1536$ discretization points per ellipsoid)
bounded by a box two wavelengths in
size takes over $1000$ GMRES iterations to get the residual below $10^{-3}$
(see Figure~\ref{figure:ellipsoid_party} and Remark \ref{re:gmres}).
Moreover, iterative methods are unable to efficiently
handle problems with multiple right hand sides that arise frequently in design
problems.

In this manuscript, we describe our contributions towards providing a
remedy to problems (2) and (3). We consider scattering problems
that arise in applications such as sonar and radar where the goal is
to find the scattered field at a set of locations $\Gamma_{\rm out}$ generated
when an incident wave sent from locations $\Gamma_{\rm in}$ hit an object or collections
of objects (whose surface we denote $\Gamma$), cf.~Fig.~\ref{fig:geom}(a).  The locations
of the source and target points are known and typically given by the user.
In this context, we present an efficient high-order accurate method for constructing
a scattering matrix $\mtx{S}$ which maps the incident field generated by the given
sources on $\Gamma_{\rm in}$ to a set of induced charges on $\Gamma$ that generate the scattered
field on $\Gamma_{\rm out}$.
A crucial component to constructing $\mtx{S}$ is solving a boundary integral equation defined on $\Gamma$.
A recently developed high order accurate discretization technique for integral
equations on surfaces \cite{Bremer-Gimbutas1,Bremer-Gimbutas2} is used.   This
discretization has proven to be effective for Laplace and low frequency scattering problems.
Since many scattering problems of interest involve higher frequencies and complicated
geometries, a large number of discretization points are typically required to achieve high accuracy,
even for high-order methods. In consequence, it would be
infeasible to solve the resulting linear systems using dense linear algebra (e.g., Gaussian elimination).
As a robust alternative to the commonly used iterative solvers (whose convergence problems in
the present context was previously discussed), we propose the use of a
direct solver that exploits the underlying structure of the matrices of interest
to construct an approximation to the scattering matrix $\mtx{S}$. Direct solvers
have several benefits over iterative solvers:
(i) the computational cost does not depend on the spectral properties of the system,
(ii) each additional incoming field can be processed quickly,
(iii) the amount of memory required to store the scattering matrix scales linearly with the
number of discretization points.

\subsection{Model problems}
\label{sec:modprob}
The solution technique presented is applicable to many boundary integral equations
that arise in mathematical physics but we limit our discussion to
the boundary value problems presented in this section.

Let $\Omega$ denote a scattering body, or a union of scattering bodies;
Figures \ref{figure:ellipsoid_party}, \ref{figure:deformed_torus}, \ref{figure:heart}, \ref{figure:cube}, and
\ref{fig:trefoil} of Section~\ref{sec:numerics} illustrate some specific domains of interest.
Let $\Gamma = \partial\Omega$ denote the surface of the body (or bodies).
We are interested in finding the scattered field $u$ off $\Omega$
for a given incident wave $u_{\rm in}$.
Specifically this involves solving either an exterior Dirichlet problem for \emph{sound-soft}
scatterers
\begin{equation}
\label{formulations:dirichlet}
\tag{DBVP}
\left\{
\begin{aligned}
\Delta u + \kappa^2 u =&\ 0\qquad&&\mbox{ in } \Omega^c,\\
                    u =&\ g\qquad&&\mbox{ on } \Gamma,\\
\frac{\partial u}{\partial r}  - i\kappa  u  =&\ O  \left( \frac{1}{r} \right)\qquad&&\mbox{ as } r\to\infty,
\end{aligned}\right.
\end{equation}
where $g = -u_{\rm in}$, or an exterior Neumann problem for \emph{sound-hard} scatterers
\begin{equation}
\label{formulations:neumann}
\tag{NBVP}
\left\{
\begin{aligned}
\Delta u + \kappa^2 u =&\ 0\qquad&&\mbox{ in } \Omega^c,\\
\frac{\partial u}{\partial \eta} =&\ g\qquad&&\mbox{ on } \Gamma,\\
\frac{\partial u}{\partial r}  - i\kappa  u  =&\ O  \left( \frac{1}{r} \right)\qquad&&\mbox{ as } r\to\infty,
\end{aligned}\right.
\end{equation}
where $\eta$ is the outwards point unit normal vector to $\Gamma$, and where
$g = -\partial_{\eta}u_{\rm in}$.

\subsection{Outline}
The paper proceeds by introducing the integral formulations for each choice of boundary condition
in Section \ref{sec:form}.  Section \ref{sec:disc} presents a high order accurate technique for
discretizing the integral equations.  Section \ref{sec:scatteringmatrices} formally defines a scattering
matrix and describes a method for constructing an approximation of the scattering matrix.
Section \ref{sec:hierarchical_solver} describes an accelerated technique for constructing the scattering matrix
that often has asymptotic cost of $O(N^{3/2})$ (providing the scattering body and the wave-length
of the radiating field are kept constant as $N$ increases).
Section \ref{sec:numerics} presents the results of several numerical experiments
which demonstrate the properties of the proposed solver.

\end{section}

\begin{section}{Integral equation formulations of linear elliptic boundary value problems}
\label{sec:form}
The reformulation of the boundary value problems in Section \ref{sec:modprob} as
integral equations involves the classical single and double layer operators
\begin{equation}
\begin{split}
S_\kappa f(x)   &=   \dblint_{\Gamma} G_\kappa (x,y) f(y) ds(y),                        \\
D_\kappa f(x)   &=   \dblint_{\Gamma} \eta_y \cdot \nabla_y G_\kappa (x,y) f(y) ds(y)\ \ \mbox{and} \\
D_\kappa ^*f(x) &=   \dblint_{\Gamma} \eta_x \cdot \nabla_x G_\kappa (x,y) f(y) ds(y)
\end{split}
\label{formulations:operators}
\end{equation}
where $G_\kappa (x,y)$ is the free space Green's function
\begin{equation*}
G_\kappa (x,y) = \frac{\exp\left(i\kappa  \left|x-y\right|\right)}{\left|x-y\right|}
\end{equation*}
associated with the Helmholtz equation at wavenumber $\kappa $
and where $\eta_p$ denotes the outward-pointing unit normal to the surface $\Gamma$
at the point $p$.

The reformulation of (\ref{formulations:dirichlet}) and (\ref{formulations:neumann})
is not entirely straightforward because of the possibility
of spurious resonances.  That is, for certain values of the wavenumber $\kappa$, the na\"{\i}ve integral equation
formulations of (\ref{formulations:dirichlet}) and (\ref{formulations:neumann}) lead to operators
with nontrivial null-spaces, even though the original boundary value problems are well-posed.
The wave-numbers $\kappa $ for which this occurs are called {\it spurious resonances}.
Complications arise due to the existence of multiple solutions at a resonance.
An additional (and perhaps more serious) problem is that for values of $\kappa$ close to resonant values,
the integral operators (\ref{formulations:operators}) become ill-conditioned.

\subsection{The boundary integral formulation for Dirichlet boundary values problems }
\label{sec:dir}
The double layer representation
\begin{equation*}
u(x) = D_\kappa  \sigma(x).
\end{equation*}
of the solution of (\ref{formulations:dirichlet}) leads to the
second kind integral equation
\begin{equation*}
\frac{1}{2}\sigma(x) + D_\kappa  \sigma(x) = g(x).
\end{equation*}
As is well known, this integral equation is not be uniquely solvable in the event that the eigenproblem
\begin{equation*}
\begin{split}
\Delta u + \kappa ^2 u = 0 \ \ \  & \mbox{ in } \Omega                         \\
\frac{\partial u}{\partial \eta} = 0 \ \ \  & \mbox{ on } \Gamma   \\
\end{split}
\end{equation*}
admits a nontrivial solution.  A uniquely solvable integral equation can be obtained by choosing to represent the
solution as
\begin{equation}
\label{eq:combine}
u(x) = D_\kappa  \sigma(x) + i|\kappa| S_\kappa  \sigma(x);
\end{equation}
see, for instance, \cite{Nedelec,Colton-Kress2}.  This leads to the
integral equation
\begin{equation}
\frac{1}{2} \sigma(x)  + D_\kappa \sigma(x) + i |\kappa|  S_\kappa  \sigma(x) = g(x),
\label{formulations:combined}
\end{equation}
which is sometimes referred to as the Combined Field Integral
Equation (CFIE).  We make use of the CFIE to solve the
problem (\ref{formulations:dirichlet}) in this article.

\subsection{The boundary integral formulation for Neumann boundary values problems }
\label{sec:neu}

%
The single layer representation
\begin{equation*}
u = S_\kappa  \sigma.
\end{equation*}
of the solution of (\ref{formulations:neumann}) leads
to the second kind integral equation
\begin{equation}
\frac{1}{2} \sigma(x) + D_\kappa ^*\sigma(x) = g(x).
\label{formulations:naive}
\end{equation}
Once again the equation (\ref{formulations:naive}) is not necessarily
uniquely solvable --- when the eigenproblem
\begin{equation}
\begin{split}
\Delta u + \kappa ^2 u = 0 \ \ \  & \mbox{ in } \Omega                         \\
u = 0 \ \ \  & \mbox{ on } \Gamma   \\
\end{split}
\label{eq:eig}
\end{equation}
has nontrivial solutions (that is, when $\kappa$ is a resonant wavenumber), the operator appearing on the right-hand side of
(\ref{formulations:naive}) is not invertible.
Since the focus of this article is on the discretization of
integral operators and the rapid inversion of the resulting
matrices, we simply avoid wavenumbers $\kappa$ for which
spurious resonances exist.  Appendix \ref{app:nbie} presents a robust
integral equation formulation which is suitable for use by the solver
of this paper.    We forgo its use here in order to reduce the complexity
of the presentation.


\end{section}

\begin{section}{A Nystr\"om method for the discretization of integral equations on surfaces}
\label{sec:disc}
In this work, we use a modified Nystr\"om method for the
discretization of the integral operators (\ref{formulations:operators}) of the proceeding section.
Our approach differs from standard Nystr\"om  methods in two important ways:

\begin{enumerate}[(1)]

\item
It captures the $L^2$ action of an operator rather than its pointwise behavior.
Among other advantages, this means that the approach results in well-conditioned
operators in the case of Lipschitz domains.

\item
A highly accurate mechanism for evaluating the singular integrals which arise
in Nystr\"om discretization whose cost is largely independent of the geometry
of the surface is employed.  This is in contrast to standard methods for the
evaluation of the singular integrals arising from the Nystr\"om discretization
of integral operators on surfaces which involve constants which depend
on the geometry of the domain.
\end{enumerate}

This manuscript gives a cursory account of the method.
A detailed description  can be found in
\cite{Bremer-Gimbutas1} and a discussion of the advantages
of $L^2$ discretization over standard Nystr\"om and collocation
methods can be found in \cite{Bremer1}.

\begin{subsection}{Decompositions and discretization quadratures.}

We define a decomposition $D$ of a surface $\Gamma \subset \mathbb{R}^3$
to be a finite sequence
\begin{equation*}
 \{ \rho_j: \Delta^1 \to \Gamma \}_{j=1}^{m}
\end{equation*}
of smooth mappings given on the simplex
\begin{equation*}
\Delta^1 = \left\{
\left(x,y\right) \in \mathbb{R}^2 \ :\
0 \leq x \leq 1, \ 0\leq y \leq 1-x
\right\}
\end{equation*}
with non-vanishing Jacobian determinants
such that the sets
\begin{equation*}
\rho_1(\Delta^1), \rho_2(\Delta^1), \ldots, \rho_m(\Delta^1)
\end{equation*}
form a disjoint union of $\Gamma$.

We call a quadrature rule
$\{x_1,\ldots,x_l,w_1,\ldots,w_l\}$
on $\Delta^1$ with positive weights which integrates all elements of
the space $\mathscr{P}_N$ of polynomials of degree less than or equal to
$2N$ on the simplex $\Delta^1$ a {\it discretization quadrature.}
Such a rule is called Gaussian if
the length $l$ of the rule is equal to $(N+1)\cdot(N+2)/2$.  That is, a Gaussian
rule is one whose length is equal to the dimension of the space of polynomials
of degree less than or equal to $N$ but which integrates polynomials
of degree less than or equal to $2N$.
For the most part, Gaussian quadratures on triangles are not available.  In the experiments
of this paper, we used quadrature rules constructed
via a generalization of the procedure of \cite{Bremer-Gimbutas-Rokhlin}.
Table~\ref{table:discretization_rules} lists the properties of the
discretization quadrature rules used in this paper and compare their lengths
to the size of a true Gaussian quadrature.


\begin{table}[h]
\begin{tabular}{lr|rrrrrr}
Discretization order                 & $N$            &  4 &  6 &  8 & 10 &  12 &  16 \\
Integration order                    & $2N$           &  8 & 12 & 16 & 20 &  24 &  32 \\
Length of a true Gaussian quadrature & $(N+1)(N+2)/2$ & 15 & 28 & 45 & 66 &  91 & 153 \\
Length of the quadrature that we use & $l$            & 17 & 32 & 52 & 82 & 112 & 192
\end{tabular}

\vspace{1.5mm}

\caption{The properties of the quadratures used in this paper. We compare the length of
an ideal (hypothetical) Gaussian quadrature to the quadrature that we actually have.}
\label{table:discretization_rules}
\end{table}

\end{subsection}

\begin{subsection}{Discretizations of integral operators.}
We now associate with a decomposition
$$D=\{\rho_j : \Delta^1 \to \Gamma\}_{j=1}^m$$
and a discretization quadrature
$\{x_1,\ldots,x_l,w_1,\ldots,w_l\}$
a scheme for the discretization of certain integral operators.
We begin by letting $S$ be the subspace of $L^2(\Gamma)$
consisting of all functions $f$ which are pointwise
defined on $\Gamma$ and such that for each $j=1,\ldots,m$
the function
\begin{equation*}
f(\rho_j(x))\left| d\rho_j(x)^*  d\rho_j(x)\right|^{1/2}
\end{equation*}
is a polynomial of order $N$ on $\Delta^1$.
Denote by $P$ be the projection of $L^2(\Gamma)$ onto the subspace $S$
and let  $\Phi: S \to \mathbb{C}^{ml}$ be a mapping which
takes $f \in S$ to a vector with entries
\begin{equation*}
f\left(\rho_j(x_i)\right)\sqrt{w_i} \left| d\rho_j(\rho_j(x_i))^*
d\rho_j(\rho_j(x_i))\right|^{1/4},
\ \ \  j=1,\ldots,m,\ \ i=1,\ldots,l.
\end{equation*}
The ordering of the entries of $\Phi(f)$ is irrelevant.
Suppose that $T:L^2(\Gamma)\to L^2(\Gamma)$ is of the form
\begin{equation}
Tf(x) =  \dblint_{\Gamma} K(x,y) f(y) ds(y),
\label{nystrom:intop}
\end{equation}
with $K$ a linear combination of the kernels listed in
(\ref{formulations:operators}).
Let $A$ denote the $ml\times ml$ matrix which commutes
the discretization of the operator $T$ induced by the specified decomposition
and discretization quadrature as specified in the diagram (\ref{discretization:commutative}).
\begin{equation}
\begin{CD}
S \subset L^2(\Gamma)@>PT >> S \subset L^2(\Gamma)\\
@VV\Phi V   @ VV\Phi V \\
\mathbb{C}^{ml} @>A>>\mathbb{C}^{ml}
\end{CD}
\label{discretization:commutative}
\end{equation}

\end{subsection}

\begin{subsection}{Quadrature.}
\label{sec:quadrature}
Forming the matrix $A$ which appears the diagram (\ref{discretization:commutative})
is an exercise in numerical quadrature.
Let $x$ be a discretization node on $\Delta$ and let $w$ be
the corresponding weight.  The $l$ entries of the matrix $A$
in (\ref{discretization:commutative}) associated with the mapping $\rho_j$ and the point $\rho_i(x)$ must map the scaled values
\begin{equation}
f(\rho_j(x_1))|d\rho^*_j(x_k)d\rho_j(x_k)|^{1/4}\sqrt{w_k},\qquad k = 1,2,\dots,l,
\label{nystrom:scaled}
\end{equation}
of the function $f$ at the images of the discretization quadrature nodes under $\rho_j$ to the
value of the integral
\begin{equation}
\dblint_{\Delta^1} K(\rho_i(x),\rho_j(y)) f(\rho_j(y)) |d\rho^*_j(y)d\rho_j(y)|^{1/2} \sqrt{w} \ dy.
\label{nystrom:int}
\end{equation}
%

The manner in which this vector is formed depends
on the location of the point $\rho_i(x)$ vis-\`a-vis the patch $\rho_j(\Delta^1)$.  Let
$B$ denote the ball of minimum radius containing $\Gamma$ and $2B$ denote the ball with the
same center as $B$ but with twice the radius.  Then we say that a point $x$ is well-separated
from a surface $\Gamma$ if $x$ is outside of $2B$, and we say that $x$ is near $\Gamma$
if $x$ is inside of the ball $2B$ but not on the surface $\Gamma$.

\begin{subsubsection}{Smooth regime.}
When the point $\rho_i(x)$ is well-separated from $\rho_j(\Delta^1)$ the kernel $K$ is smooth and
the discretization quadrature can be used to evaluate the integral (\ref{nystrom:int}).
In this case, we take the $l$ entries of the matrix $A$ to be
\begin{equation}
K(\rho_i(x),\rho_j(x_k))\left| d\rho_j^*(x_k)d\rho_j(x_k) \right|^{1/4} \sqrt{w} \sqrt{w_k},\ \ \ k=1,2,\ldots,l.
\label{discretization:far}
\end{equation}
\end{subsubsection}

\begin{subsubsection}{Near regime.}
When $\rho_i(x)$ is near $\rho_j(\Delta^1)$,
formula (\ref{discretization:far}) may not
achieve sufficient accuracy.  Thus, we adaptively form a quadrature
with nodes and weights
\begin{equation*}
u_1,\ldots,u_m,r_1,\ldots,r_m
\end{equation*}
sufficient to evaluate the nearly singular integrals which arise
and apply an interpolation matrix $I$ which takes the scaled values
(\ref{nystrom:scaled}) of the function $f$ at the nodes
of the discretization quadrature to its scaled values at the
nodes of the adaptive quadrature.
The vector of entries of $A$ is given by the product of the vector
of values of the kernel evaluated at the adaptive quadrature nodes
with the matrix $I$.  As in the scheme of \cite{Bremer-Gimbutas1}, the adaptive quadrature
is formed in such a way that the matrix $I$  has special structure --- it is a product of block diagonal matrices ---
which allows it to be applied efficiently.
\end{subsubsection}

\begin{subsubsection}{Singular regime.}
When the target node lies on the surface $\rho_j(\Delta^1)$, the integral
(\ref{nystrom:int}) is singular and its evaluation becomes quite difficult.
The typical approach to the evaluation of integrals of this form is to apply
the Duffy transform or polar coordinate transform in order to remove
the singularity of the kernel.  Once this has been done, quadratures
for smooth functions (for instance, product Legendre quadratures) are typically applied.

While the such schemes are exponentially convergent, they can be
extremely inefficient.  The difficulty is that the functions resulting
from the change of variables can have poles close to the integration domain.
For instance, if a polar coordinate transform is applied, the singular integrals
which must be evaluated take the form
\begin{equation}
\int_0^{2\pi}
\hskip .1em
\int_0^{{R}(\theta)}
\left(q_{-1}(\theta) + q_{0}(\theta)r + q_{1}(\theta) r^2   +
q_{2}(\theta) r^3   + \cdots\ \right)
dr\hskip .1em
d\theta
\label{singular:polar2}
\end{equation}
where the function $R(\theta)$ gives a parameterization of the integration
domain $\Delta^1$ in polar coordinates and each $q_j(\theta)$ is a
quotient
form
\begin{equation*}
\frac{r_j(\theta)}{l(\theta)^{\alpha}},
\end{equation*}
where $\alpha$ is a positive constant dependent on the order $j$
and the kernel under consideration, $r_j(\theta)$ is a trigonometric polynomial of finite order,
and $l(\theta)$ is a function which depends on the parameterization $\rho_j$.
More specifically, if
\begin{equation*}
\xi_1 = \frac{\partial \rho_j (x_1,x_2)}{\partial x_1} \ \ \ \mbox{ and }\ \ \
\xi_2 = \frac{\partial \rho_j(x_1,x_2)}{\partial x_2},
\end{equation*}
then  the function $l(\theta)$ is given by
\begin{equation*}
l(\theta) = \sqrt{ \lambda \cos^2(\theta) + 2 \gamma \cos(\theta) \sin(\theta) + \lambda^{-1} \sin^2(\theta) }
\end{equation*}
where
\begin{equation*}
\cos(\gamma) = \frac{\xi_1 \cdot \xi_2 }{\left|\xi_1\right| \left|\xi_2\right|}
\ \ \ \mbox{ and }\ \ \
\lambda =
\frac{\left|\xi_1\right|}{\left|\xi_2\right|}.
\end{equation*}
Obviously, the function $l(\theta)$ can have zeros close to the real axis.

The approach used in this paper --- and described in \cite{Bremer-Gimbutas1} --- calls
for first applying a change of variables in order to ensure that the mapping $\rho_j$
is conformal at the target node.  Then the function $l(\theta)$ is a constant
which does not depend on $\theta$ and the integrand simplifies considerably.  The unfortunate
side effect of applying this change of variables is that the integration domain is no
longer the standard simplex $\Delta^1$ but rather an arbitrary triangle.  Now
the integral which must be evaluated takes the form
\begin{equation}
\int_0^{2\pi}
\hskip .1em
\int_0^{\tilde{R}(\theta)} \left(
\tilde{q}_{-1}(\theta) + \tilde{q}_{0}(\theta)r + \tilde{q}_{1}(\theta) r^2   +
\tilde{q}_{2}(\theta) r^3   + \cdots \ \right)
dr\hskip .1em
d\theta
\label{singular:polar3}
\end{equation}
with the $\tilde{q_j}$ trigonometric polynomials of finite order.
The function $\tilde{R}(\theta)$ can have poles close to the integration domain, however.
We combat this problem by precomputing a collection of quadrature rules
designed to integrate smooth functions on arbitrary triangles and apply them to
evaluate the integrals (\ref{singular:polar3}).  The performance of this approach
is described in \cite{Bremer-Gimbutas1}; in one example of that paper, the standard
approach of applying a change of variables followed by Gauss-Legendre quadrature
required more than $100,000$ quadrature nodes to accurately evaluate the singular
integrals arising from a boundary value problem while the precomputed rules used
here and in that paper required less than $1,000$ nodes.
\end{subsubsection}

\end{subsection}

\label{section:nystrom}

\end{section}

\begin{section}{Scattering matrices as a numerical tool}
\label{sec:scatteringmatrices}
This section describes the concept of a scattering matrix, first from a physical
perspective, and then the corresponding linear algebraic interpretation. It also
informally describes how to build approximate scattering matrices computationally.
A rigorous description is given in Section \ref{sec:hierarchical_solver}.

\subsection{The scattering operator}
\label{sec:physical}
We consider a ``sound-soft'' scattering formulation illustrated in Figure \ref{fig:geom}.
We are given a charge distribution $s_{\rm in}$ on a contour $\Gamma_{\rm in}$
(the radiation ``source'') that generates an incoming field $g$ on the scatterer $\Gamma$.
The incoming field on the scatterer induces an outgoing field $u$ that satisfies
(\ref{formulations:dirichlet}).
Our objective is to computationally construct $u_{\rm out}$, the restriction of the outgoing field
to the specified contour $\Gamma_{\rm out}$, which represents the ``antenna.''
To solve (\ref{formulations:dirichlet}), we use the combined field formulation (\ref{eq:combine}).
In other words, we represent the outgoing field via a distribution $\sigma$ of ``induced sources''
on the contour $\Gamma$. The map we seek to evaluate can be expressed mathematically as
\begin{equation}
\label{eq:3step_map}
\begin{array}{ccccccccc}
u_{\rm out} &=&
T_{\rm out} & T_{\Gamma}^{-1} & T_{\rm in} & s_{\rm in}, \\[2mm]
\mbox{\textit{outgoing field on }}\Gamma_{\rm out}
&&&&&
\mbox{\textit{given charges on }}\Gamma_{\rm in}
\end{array}
\end{equation}
where $T_{\rm in}\,\colon\,L^{2}(\Gamma_{\rm in}) \rightarrow L^{2}(\Gamma)$ is the operator
$$
[T_{\rm in}s_{\rm in}](x) = \int_{\Gamma_{\rm in}} G_{\kappa}(x,y) s_{\rm in}(y) ds(y),\qquad x \in \Gamma,
$$
 $T_{\Gamma}\,\colon\,L^{2}(\Gamma) \rightarrow L^{2}(\Gamma)$ is the operator
\begin{equation}
\label{eq:T_gamma}
[T_{\Gamma}g](x) =
\frac{1}{2} g(x)  + D_\kappa g(x) + i |\kappa| S_\kappa g(x),\qquad x \in \Gamma,
\end{equation}
and $T_{\rm out}\,\colon\,L^{2}(\Gamma)\rightarrow L^{2}(\Gamma_{\rm out})$ is the operator
$$
[T_{\rm out}\sigma](x) = \int_{\Gamma} G_\kappa(x,y)\sigma(y)ds(y) +\int_{\partial\Omega} \partial_{\nu_y}G_\kappa(x,y)\sigma(y)ds(y),
\qquad x \in \Gamma_{\rm out}.
$$
Observe that since $T_{\Gamma}$ is an invertible second kind integral operator, its inverse
has singular values bounded away from zero.

Given a finite precision $\varepsilon > 0$,
we say that an operator $S\,\colon\,L^{2}(\Gamma) \rightarrow L^{2}(\Gamma)$ is a ``scattering operator''
for $\Gamma$ to within precision $\varepsilon$ when
$$
\|T_{\rm out}\,T_{\Gamma}^{-1}\,T_{\rm in} - T_{\rm out}\,S\,T_{\rm in}\| \leq \varepsilon.
$$
It turns out that the operators $T_{\rm in}$ and $T_{\rm out}$ have rapidly decaying singular values.
Let $P_{\rm in}^{(k)}$ denote the orthogonal projection onto
the set of $k$ leading left singular vectors of $T_{\rm in}$, and let $P_{\rm out}^{(k)}$ denote the orthogonal
projection onto the set of $k$ leading right singular vectors of $T_{\rm out}$. Then
$$
\|T_{\rm out}P_{\rm out}^{(k)} - T_{\rm out}\| \leq \sigma_{k+1}(T_{\rm out})
\qquad\mbox{and}\qquad
\|T_{\rm  in}P_{\rm  in}^{(k)} - T_{\rm  in}\| \leq \sigma_{k+1}(T_{\rm in }).
$$
It follows that if we define an approximate scattering operator $S$ via
\begin{equation}
\label{eq:def_S_cont}
S = P_{\rm out}^{(k)}\,T_{\Gamma}^{-1}\,P_{\rm in}^{(k)},
\end{equation}
then the approximation error
$$
\|T_{\rm out}\,T_{\Gamma}^{-1}\,T_{\rm in} - T_{\rm out}\,P_{\rm out}^{(k)}\,T_{\Gamma}^{-1}\,P_{\rm in}^{(k)}\,T_{\rm in}\|
$$
converges to zero very rapidly as $k$ increases.

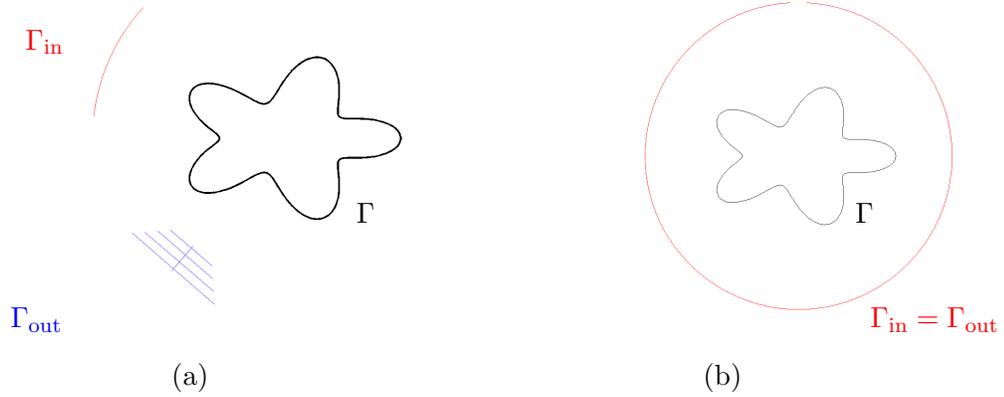
\begin{figure}
\begin{tabular}{ccc}
\setlength{\unitlength}{1mm}
\begin{picture}(55,50)
\put(00,00){\includegraphics[height=50mm]{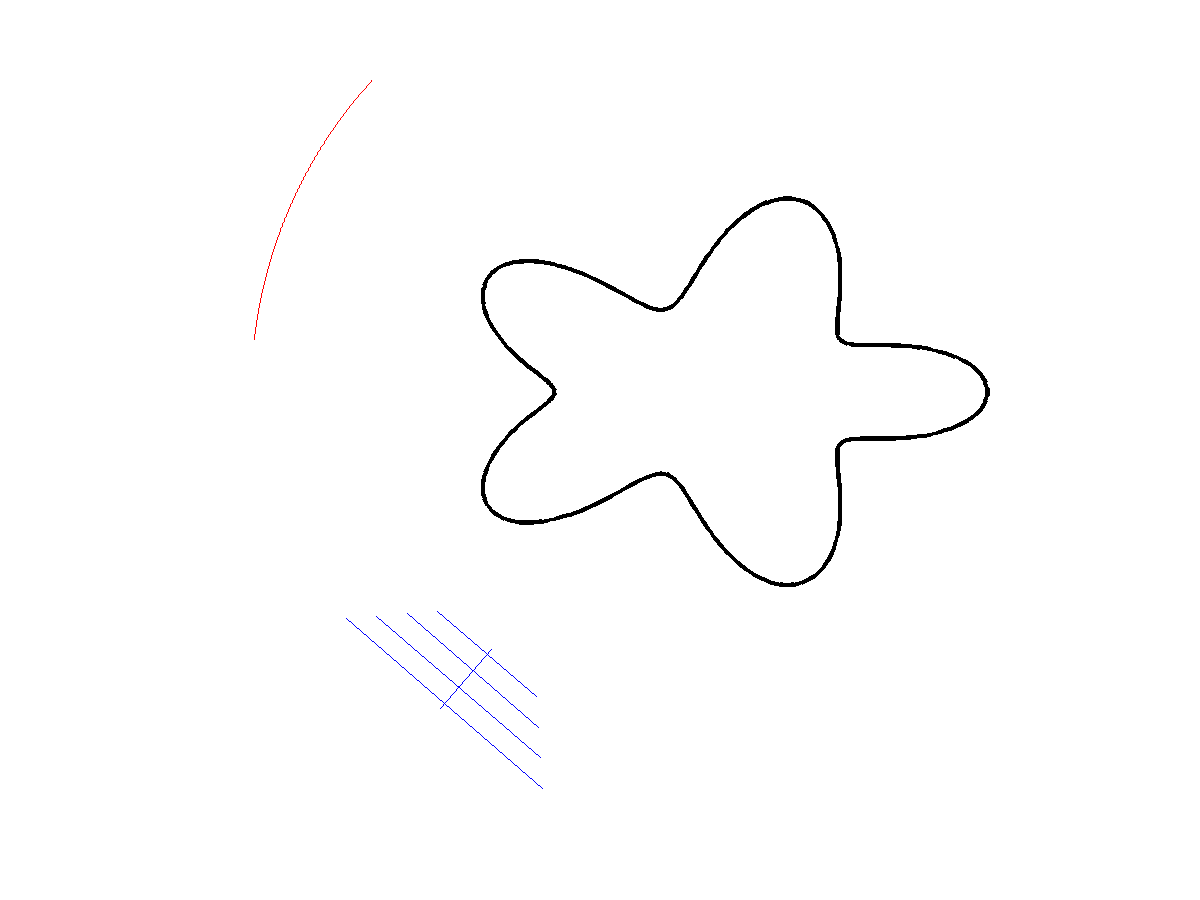}}
\put(03,03){\color{blue}$\Gamma_{\rm out}$}
\put(05,40){\color{red}$\Gamma_{\rm in}$}
\put(49,17){$\Gamma$}
\end{picture}
&\mbox{}\hspace{10mm}\mbox{}&
\setlength{\unitlength}{1mm}
\begin{picture}(50,50)
\put(00,00){\includegraphics[height=50mm]{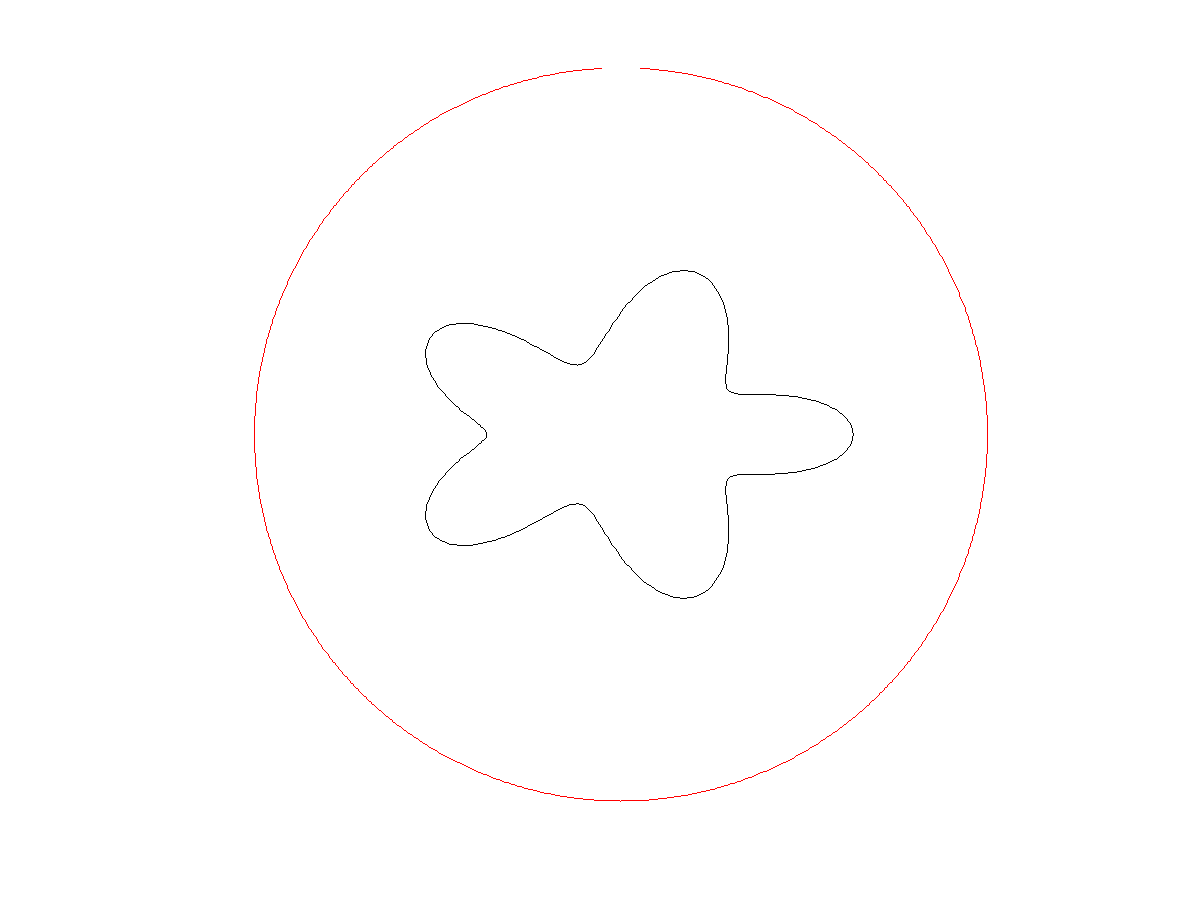}}
\put(42,17){$\Gamma$}
\put(44,03){\color{red}$\Gamma_{\rm in} = \Gamma_{\rm out}$}
\end{picture}\\
(a) && (b)
\end{tabular}
\caption{Geometry of the scattering problem in Section \ref{sec:scatteringmatrices}.
The scattering surface is $\Gamma$.
(a) The radiation source $\Gamma_{\rm in}$ is different from the antenna $\Gamma_{\rm out}$.
(b) The typical geometry considered with $\Gamma_{\rm in} = \Gamma_{\rm out}$. The scattering
matrix generated by this geometry is the ``complete'' scattering matrix in the sense that any
incoming field can be handled (as long as the radiation source is not inside $\Gamma_{\rm in}$),
and the full radiated field can be obtained.}
\label{fig:geom}
\end{figure}


\subsection{Definition of a discrete scattering matrix}
\label{sec:disc_scatter_matrix}
The continuum scattering problem described in Section \ref{sec:physical}
has a direct linear algebraic analog for the discretized equations.
Forming the Nystr\"om discretization of the three operators in
(\ref{eq:3step_map}), we obtain the map
\begin{equation}
\label{eq:3step_map_disc}
\begin{array}{ccccccccc}
\uu_{\rm out} &=&
\mtx{A}_{\rm out} & \mtx{A}^{-1} & \mtx{A}_{\rm in} & \vct{s}_{\rm in}. \\[2mm]
N_{\rm out} \times 1 && N_{\rm out} \times N & N\times N & N\times N_{\rm in} & N_{\rm in} \times 1
\end{array}
\end{equation}
Suppose that $\mtx{A}$ is an $N\times N$ matrix, and let $k_{1}$ denote a bound on the $\varepsilon$-ranks of
$\mtx{A}_{\rm out}$ and $\mtx{A}_{\rm in}$. Then form matrices $\mhU_{1}$ and $\mhV_{1}$ (the reason for the
subscripts will become clear shortly) of size $N\times k_{1}$ such that
$$
\|\mtx{A}_{\rm in } - \mhU_{1}\,\mhU_{1}^{\dagger}\,\mtx{A}_{\rm in }\| \leq \varepsilon,
\qquad\mbox{\and}\qquad
\|\mtx{A}_{\rm out} - \mtx{A}_{\rm out}\,(\mhV_{1}\mhV_{1}^{\dagger})^{*}\| \leq \varepsilon
$$
where a matrix with the superscript $\dagger$ denotes the pseudo-inverse of the matrix.  In other words, the columns
of $\mhU_{1}$ span the column space of $\mtx{A}_{\rm in}$ and the
columns of $\mhV_{1}$ span the row space of $\mtx{A}_{\rm out}$. Then the map
(\ref{eq:3step_map_disc}) can be approximated as
\begin{equation}
\begin{array}{ccccccccc}
\uu_{\rm out} &=&
\bigl(\mtx{A}_{\rm out}\,(\mhV_{1}^{\dagger})^{*}\bigr) & \mtx{S}_{1} &
\bigl(\mhU_{1}^{\dagger}\,\mtx{A}_{\rm in}\bigr) & \vct{s}_{\rm in}, \\[2mm]
N_{\rm out} \times 1 && N_{\rm out} \times k_{1} & k_{1}\times k_{1} & k_{1}\times N_{\rm in}
\end{array}
\end{equation}
where
\begin{equation}
\label{eq:def_S1}
\mtx{S}_{1} = \mhV_{1}^{*}\,\mtx{A}^{-1}\,\mhU_{1}.
\end{equation}
The $k_{1}\times k_{1}$ matrix $\mtx{S}_{1}$ is the discrete analog of the continuum scattering matrix $S$
for $\Gamma$ defined by (\ref{eq:def_S_cont}).

\subsection{Hierarchical construction of discrete scattering matrices}
\label{sec:hierarchical_intuition}
Observe that direct evaluation of the definition (\ref{eq:def_S1}) of $\mtx{S}_{1}$
via, e.g., Gaussian elimination, would be very costly since $\mtx{A}$ is dense, and
$N$ could potentially be large (much larger than $k_1$).
To avoid this expense, we will build $\mtx{S}_{1}$ via an accelerated hierarchical
procedure. As a first step, we let the entire domain
$\Gamma = \Gamma_{1}$ denote the root of the tree, and split $\Gamma_{1}$ into two
disjoint halves, cf.~Figure \ref{fig:part}(b),
\begin{equation}
\label{eq:G1=G2cupG3}
\Gamma_{1} = \Gamma_{2} \cup \Gamma_{3}.
\end{equation}
Let $n_{2}$ and $n_{3}$ denote the number of collocation points in $\Gamma_{2}$ and $\Gamma_{3}$, respectively.
The idea is now to form scattering matrices $\mtx{S}_{2}$ and $\mtx{S}_{3}$ for
the two halves, and then form $\mtx{S}_{1}$ by ``merging'' these smaller scattering matrices.
To formalize, split the matrix $\mtx{A}$ into blocks to conform with the partitioning
(\ref{eq:G1=G2cupG3}),
$$
\mA = \mtwo{\mA_{2,2}}{\mA_{2,3}}{\mA_{3,2}}{\mA_{3,3}}.
$$
The off-diagonal blocks $\mA_{2,3}$ and $\mA_{3,2}$ have very rapidly decaying
singular values,
 which means that they can be factored
\begin{align}
\label{eq:factor23}
&\begin{array}{cccccc}
\mA_{2,3} &=& \mhU_{2} & \mtA_{2,3} & \mhV_{3}^{*},\\
n_{2} \times n_{3} && n_{2} \times k_{2} & k_{2}\times k_{3} & k_{3}\times n_{3}
\end{array}
\\
\label{eq:factor32}
&\begin{array}{cccccc}
\mA_{3,2} &=& \mhU_{3} & \mtA_{3,2} & \mhV_{2}^{*}.\\
n_{3} \times n_{2} && n_{3} \times k_{3} & k_{3}\times k_{2} & k_{2}\times n_{2}
\end{array}
\end{align}
The matrices $\{\mhU_{2},\,\mhU_{3}\}$ and $\{\mhV_{2},\,\mhV_{3}\}$ are
additionally assumed to span the column spaces of $\mtx{A}_{\rm in}$ and $\mtx{A}_{\rm out}^{*}$,
in the sense that there exist ``small'' matrices $\mU_{1}$ and $\mV_{1}$ such that
\begin{equation}
\label{eq:hbases}
\begin{array}{cccccc}
\mhU_{1} &=& \mtwo{\mhU_{2}}{\mtx{0}}{\mtx{0}}{\mhU_{3}} & \mU_{1},\\
N\times k_{1} && N\times (k_{2}+k_{3}) & (k_{2}+k_{3}) \times k_{1}
\end{array}
\quad\mbox{and}\quad
\begin{array}{cccccc}
\mhV_{1} &=& \mtwo{\mhV_{2}}{\mtx{0}}{\mtx{0}}{\mhV_{3}} & \mV_{1}.\\
N\times k_{1} && N\times (k_{2}+k_{3}) & (k_{2}+k_{3}) \times k_{1}
\end{array}
\end{equation}
Inserting (\ref{eq:factor23}), (\ref{eq:factor32}), and (\ref{eq:hbases}) into (\ref{eq:def_S1}), we find
\begin{equation}
\label{eq:S1_tele}
\mS_{1} =
\mV_{1}^{*}
\mtwo{\mhV_{2}^{*}}{\mtx{0}}{\mtx{0}}{\mhV_{3}^{*}}
\left[\begin{array}{cc}
\mtx{A}_{2,2} & \mhU_{2}\mtA_{2,3}\mhV_{3}^{*} \\
\mhU_{3}\mtA_{3,2}\mhV_{2}^{*} & \mtx{A}_{3,3}
\end{array}\right]^{-1}
\mtwo{\mhU_{2}}{\mtx{0}}{\mtx{0}}{\mhU_{3}}
\mU_{2}.
\end{equation}
Using the Woodbury formula for matrix inversion (see Lemma \ref{lemma:mergelemma}), it can be shown that
\begin{equation}
\label{eq:S1_woodbury}
\mtwo{\mhV_{2}^{*}}{\mtx{0}}{\mtx{0}}{\mhV_{3}^{*}}
\left[\begin{array}{cc}
\mtx{A}_{2,2} & \mhU_{2}\mtA_{2,3}\mhV_{3}^{*} \\
\mhU_{3}\mtA_{3,2}\mhV_{2}^{*} & \mtx{A}_{3,3}
\end{array}\right]^{-1}
\mtwo{\mhU_{2}}{\mtx{0}}{\mtx{0}}{\mhU_{3}}
=
\left[\begin{array}{cc}
\mtx{I} & \mS_{2}\mtA_{2,3} \\
\mS_{3}\mtA_{3,2} & \mtx{I}
\end{array}\right]^{-1}\,
\left[\begin{array}{cc}
\mtx{S}_{2} & \mtx{0} \\
\mtx{0} & \mtx{S}_{3}
\end{array}\right]
\end{equation}
where $\mtx{S}_{2}$ and $\mtx{S}_{3}$ are the scattering matrices for $\Gamma_{2}$ and $\Gamma_{3}$, respectively.
The matrices $\mtx{S}_{2}$ and $\mtx{S}_{3}$ are given by the following formulas
\begin{equation}
\label{eq:S2S3}
\mtx{S}_{2} = \mhV_{2}^{*}\mA_{2,2}^{-1}\mhU_{2},
\qquad\mbox{and}\qquad
\mtx{S}_{3} = \mhV_{3}^{*}\mA_{3,3}^{-1}\mhU_{3}.
\end{equation}
Combining (\ref{eq:S1_tele}) and (\ref{eq:S1_woodbury}), we obtain the formula
\begin{equation}
\label{eq:S1_merge}
\mS_{1} =
\mV_{1}^{*}
\left[\begin{array}{cc}
\mtx{I} & \mS_{2}\mtA_{2,3} \\
\mS_{3}\mtA_{3,2} & \mtx{I}
\end{array}\right]^{-1}\,
\left[\begin{array}{cc}
\mtx{S}_{2} & \mtx{0} \\
\mtx{0} & \mtx{S}_{3}
\end{array}\right]
\mU_{1}.
\end{equation}
The upside here is that in order to evaluate (\ref{eq:S1_merge}), we only need to invert a matrix
of size $(k_{2}+k_{3})\times (k_{2}+k_{3})$, whereas (\ref{eq:def_S1}) requires inversion of a
matrix of size $N\times N$. The downside is that we are still left with the task of inverting
$\mA_{2,2}$ and $\mA_{3,3}$ to evaluate $\mtx{S}_{2}$ and $\mtx{S}_{3}$ via (\ref{eq:S2S3}).
Now, to reduce the cost of constructing $\mS_{2}$ and $\mS_{3}$, we further sub-divide
$$
\Gamma_{2} = \Gamma_{4} \cup \Gamma_{5},
\qquad\mbox{and}\qquad
\Gamma_{3} = \Gamma_{6} \cup \Gamma_{7},
$$
as shown in Figure \ref{fig:part}(c). Then it can be shown (see Lemma \ref{lemma:mergelemma}) that $\mS_{2}$
and $\mS_{3}$ are given by the formulas
$$
\mS_{2} =
\mV_{2}^{*}
\left[\begin{array}{cc}
\mtx{I} & \mS_{4}\mtA_{4,5} \\
\mS_{5}\mtA_{5,4} & \mtx{I}
\end{array}\right]^{-1}\,
\left[\begin{array}{cc}
\mtx{S}_{4} & \mtx{0} \\
\mtx{0} & \mtx{S}_{5}
\end{array}\right]
\mU_{2},
\quad\mbox{and}\quad
\mS_{3} =
\mV_{3}^{*}
\left[\begin{array}{cc}
\mtx{I} & \mS_{6}\mtA_{6,7} \\
\mS_{7}\mtA_{7,6} & \mtx{I}
\end{array}\right]^{-1}\,
\left[\begin{array}{cc}
\mtx{S}_{6} & \mtx{0} \\
\mtx{0} & \mtx{S}_{7}
\end{array}\right]
\mU_{3},
$$
where
$$
\mtx{S}_{4} = \mhV_{4}^{*}\mA_{4,4}^{-1}\mhU_{4},
\qquad
\mtx{S}_{5} = \mhV_{5}^{*}\mA_{5,5}^{-1}\mhU_{5},
\qquad
\mtx{S}_{6} = \mhV_{6}^{*}\mA_{6,6}^{-1}\mhU_{6},
\qquad
\mtx{S}_{7} = \mhV_{7}^{*}\mA_{7,7}^{-1}\mhU_{7}.
$$

\begin{figure}
\begin{tabular}{ccccc}
\setlength{\unitlength}{1mm}
\begin{picture}(40,31)
\put(00,00){\includegraphics[width=40mm]{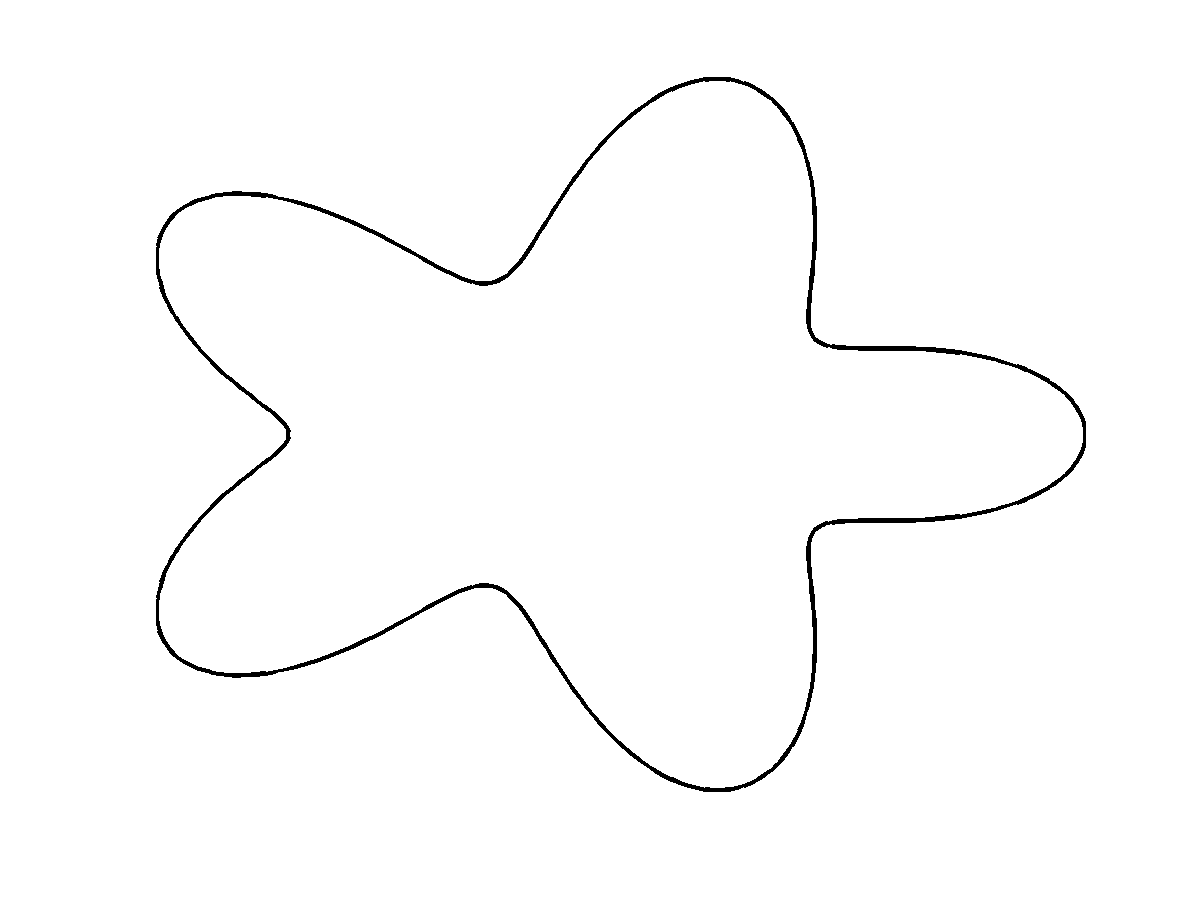}}
\put(30,00){$\Gamma = \Gamma_{1}$}
\end{picture}
&\mbox{}\hspace{10mm}\mbox{}&
\setlength{\unitlength}{1mm}
\begin{picture}(40,31)
\put(00,00){\includegraphics[width=40mm]{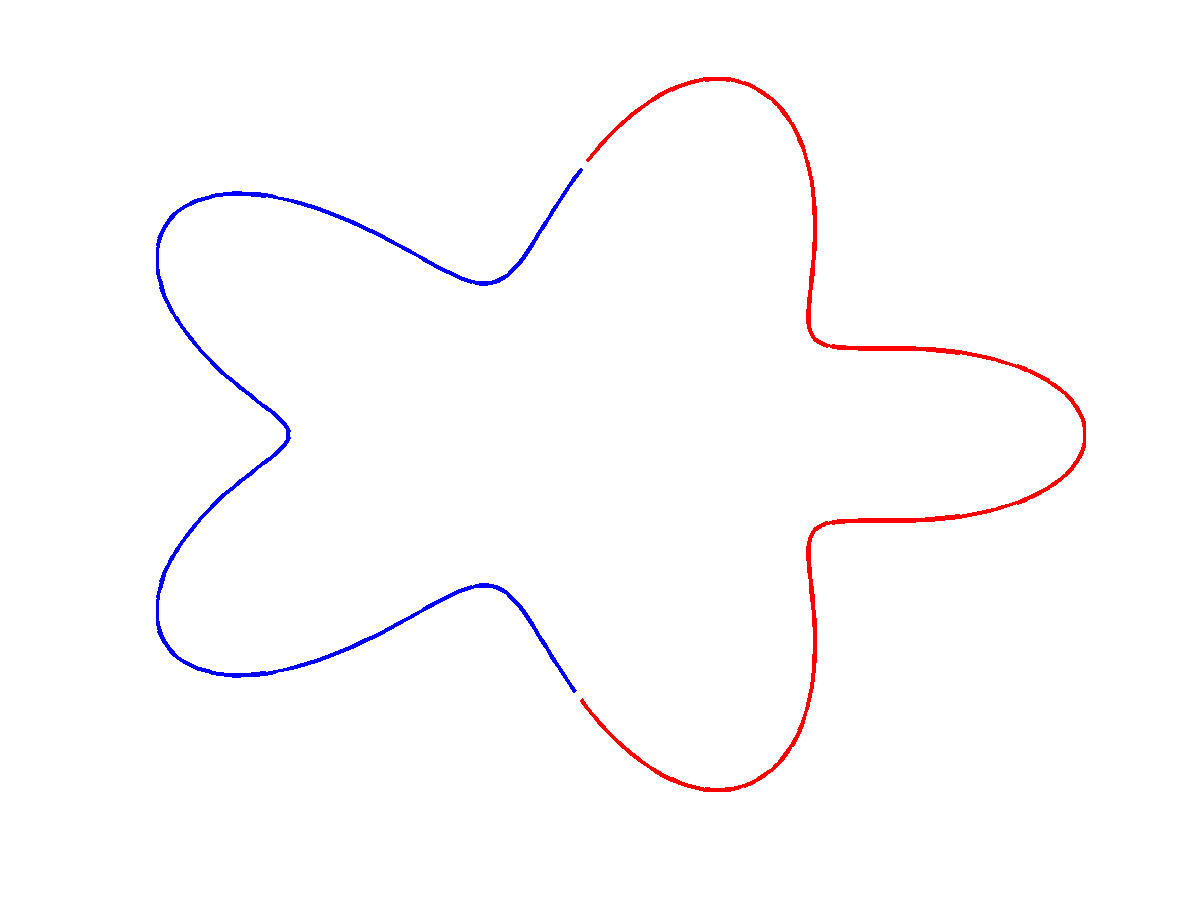}}
\put(30,00){\color{red}$\Gamma_{2}$}
\put(00,00){\color{blue}$\Gamma_{3}$}
\end{picture}
&\mbox{}\hspace{10mm}\mbox{}&
\setlength{\unitlength}{1mm}
\begin{picture}(40,31)
\put(00,00){\includegraphics[width=40mm]{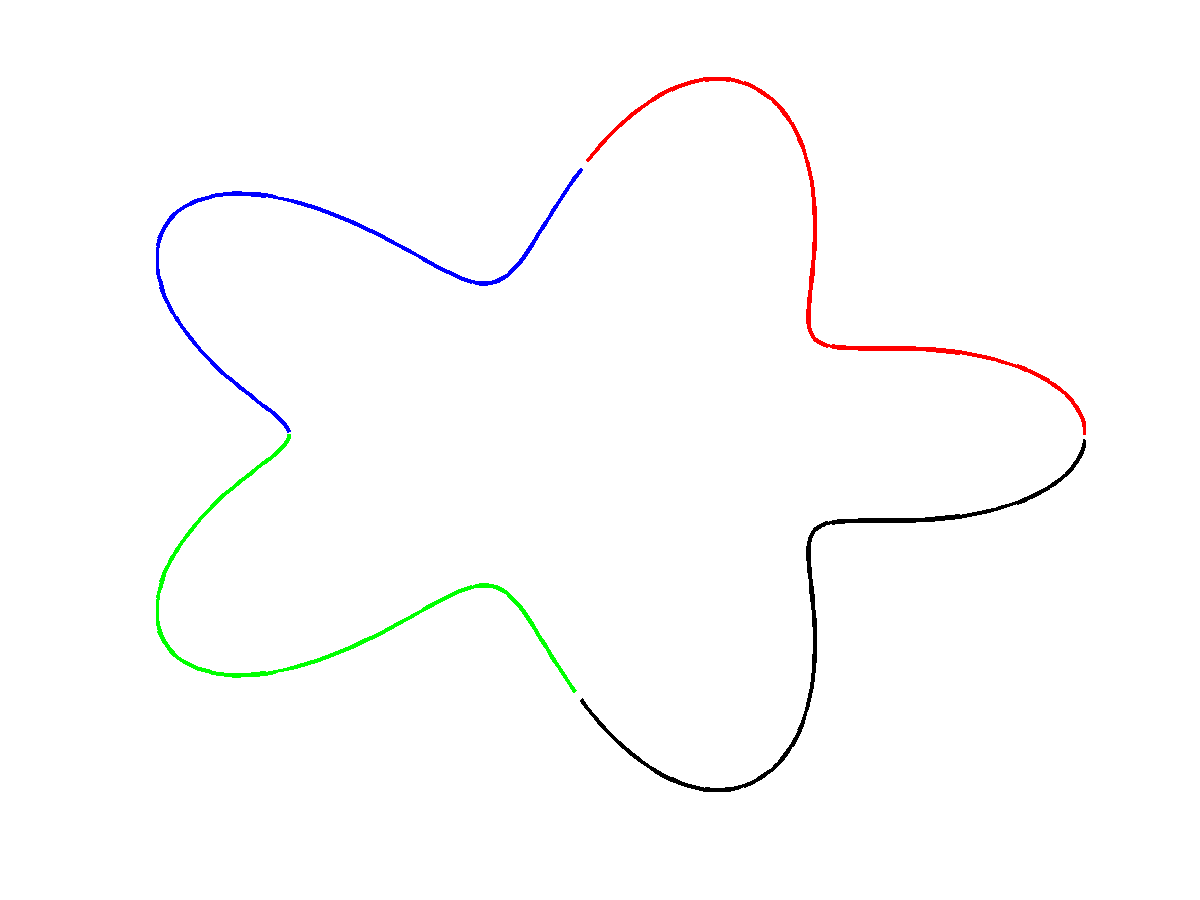}}
\put(30,00){\color{black}$\Gamma_{4}$}
\put(30,25){\color{red}  $\Gamma_{5}$}
\put(00,28){\color{blue} $\Gamma_{6}$}
\put(00,00){\color{green}$\Gamma_{7}$}
\end{picture}\\
(a) && (b) && (c)
\end{tabular}
\caption{Hierarchical partitioning of the domain as described in Section \ref{sec:hierarchical_intuition}.
(a) The full domain $\Gamma = \Gamma_{1}$.
(b) Partition $\Gamma_{1} = \Gamma_{2} \cup \Gamma_{3}$.
(c) Partition $\Gamma_{2} = \Gamma_{4} \cup \Gamma_{5}$ and $\Gamma_{3} = \Gamma_{6} \cup \Gamma_{7}$.}
\label{fig:part}
\end{figure}
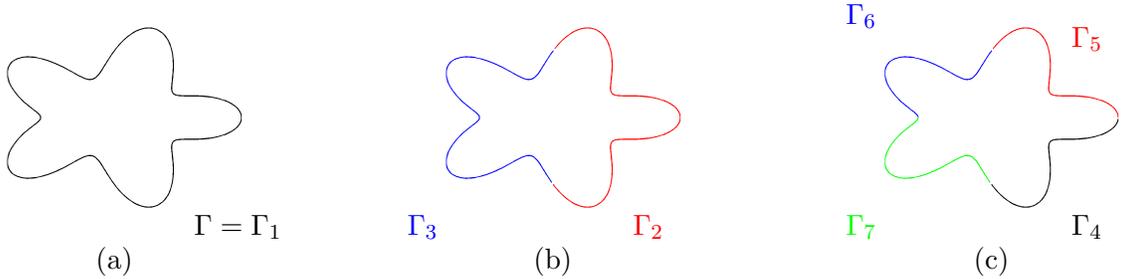

The idea is now to continue to recursively split all patches until we get to a point where each patch holds
sufficiently few discretization points that its scattering matrix can inexpensively be constructed
via, e.g., Gaussian elimination.

%

\subsection{Scattering ranks}
\label{sec:ranks}

\end{section}

\begin{section}{The hierarchical construction of a scattering matrix}
\label{sec:hierarchical_solver}
Section \ref{sec:scatteringmatrices} defined the concept of a scattering matrix,
and informally outlined a procedure for how to build approximate scattering matrices
via a hierarchical procedure. This section provides a more rigorous description, which
 requires the introduction of some notational machinery. The key concept is that the matrix
$\mtx{A}$ corresponding to the discretized operator $T_{\Gamma}$ (defined by (\ref{eq:T_gamma}))
can be represented efficiently in a data sparse format that we call \textit{Hierarchically
Block Separable (HBS)}, and is closely related to the \textit{Hierarchically Semi-Separable (HSS)}
format \cite{2009_xia_superfast,2007_shiv_sheng,2010_xia}. Our presentation is in principle
self-contained but assumes some familiarity with hierarchically rank-structured matrix formats
such as HSS/HBS (or the related $\mathcal{H}^{2}$-matrix format \cite{2010_borm_book,2004_borm_hackbusch}).

\begin{remark}
The method presented in Sections \ref{sec:disc_scatter_matrix}, \ref{sec:hierarchical_intuition}, and
this section use the same rank for the incoming and outgoing factorizations purely for simplicity
of presentation.  In practice this is not required.
\end{remark}

\subsection{Problem formulation}
Let $\mA$ denote the $N\times N$ dense matrix obtained upon Nystr\"om discretization
(see Section \ref{sec:disc}) of the operator $T_{\Gamma}$ as defined by (\ref{eq:T_gamma}).
Moreover, suppose that we are given a positive tolerance $\varepsilon$, and
two tall thin matrices $\mhU_{1}$ and $\mhV_{1}$ satisfying the following:

\vspace{1mm}

\begin{itemize}
\item[$\mhU_{1}$:] A well-conditioned matrix of size $N\times k_{1}$ whose columns
span (to within precision $\varepsilon)$ the column space of the matrix $\mtx{A}_{\rm in}$
that maps the given charge distribution on the ``radiation source'' $\Gamma_{\rm in}$ to the
collocation points on $\Gamma$. In other words, any incoming field hitting $\Gamma$ can be
expressed as a linear combination of the columns of $\mhU_{1}$.

\vspace{1mm}

\item[$\mhV_{1}$:] A well-conditioned matrix of size $N\times k_{1}$ whose columns
span (to within precision $\varepsilon)$ the row space of the matrix $\mtx{A}_{\rm out}$
that maps the induced charges on the collocation points on $\Gamma$ to the collocation
points on the ``antenna'' $\Gamma_{\rm out}$. In other words, any field generated on
$\Gamma_{\rm out}$ by induced charges on $\Gamma$ can be replicated by a charge distribution
in the span of $\mhV_{1}$.
\end{itemize}

\vspace{1mm}

Our objective is then to construct a highly accurate approximation to the scattering matrix
\begin{equation}
\label{eq:def_S1_repeat}
\mtx{S}_{1} = \mhV_{1}^{*}\,\mtx{A}^{-1}\,\mhU_{1}.
\end{equation}

\subsection{Hierarchical tree}
\label{sec:tree}
The hierarchical construction of the scattering matrix $\mS_{1}$ is
based on a binary tree structure of patches on the scattering surface $\Gamma$.
We let the full domain $\Gamma$ form the root of the
tree, and give it the index $1$, $\Gamma_{1} = \Gamma$.
We next split the root into two roughly equi-sized patches
$\Gamma_{2}$ and $\Gamma_{3}$ so that $\Gamma_{1} = \Gamma_{2} \cup \Gamma_{3}$.
The full tree is then formed by continuing to subdivide any patch
that holds more than some preset fixed number of collocation points.
A \textit{leaf} is a node in the tree corresponding to a patch that never got split.
For a non-leaf node $\tau$, its \textit{children} are the two nodes
$\alpha$ and $\beta$ such that $\Gamma_{\tau} = \Gamma_{\alpha} \cup \Gamma_{\beta}$,
and $\tau$ is then the \textit{parent} of $\alpha$ and $\beta$.
Two boxes with the same parent are called \textit{siblings}.

Let $I = I_{1} = [1,\,2,\,3,\,\dots,\,N]$ denote the full list of indices
for the collocation points on $\Gamma$, and let for any node $\tau$ the
index vector $I_{\tau} \subseteq I$ mark the set of collocation nodes inside $\Gamma_{\tau}$.
Let $n_{\tau}$ denote the number of nodes in $I_{\tau}$.

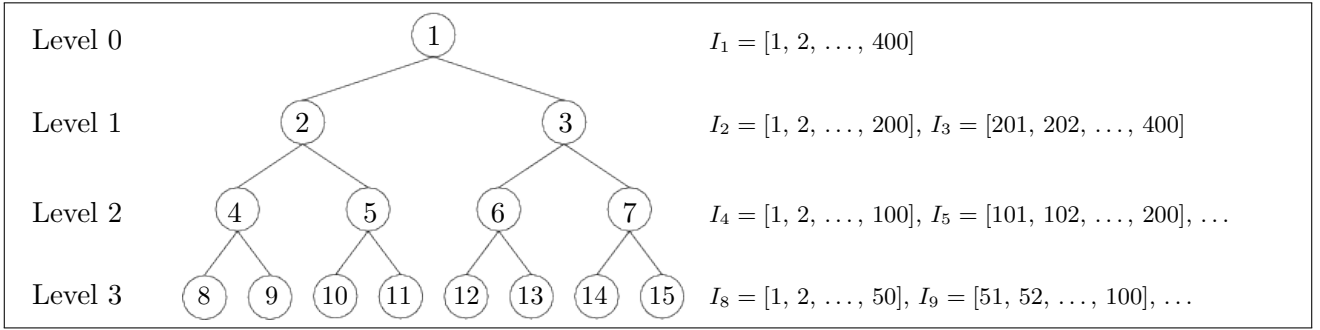
\begin{figure}
\fbox{
\setlength{\unitlength}{1mm}
\begin{picture}(169,41)
\put(20, 0){\includegraphics[height=41mm]{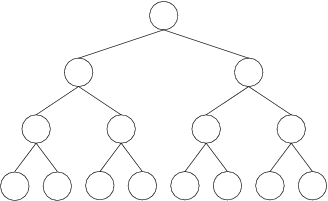}}
\put( 0,36){Level $0$}
\put( 0,25){Level $1$}
\put( 0,13){Level $2$}
\put( 0, 2){Level $3$}
\put(90,36){\footnotesize$I_{1} = [1,\,2,\,\dots,\,400]$}
\put(90,25){\footnotesize$I_{2} = [1,\,2,\,\dots,\,200]$, $I_{3} = [201,\,202,\,\dots,\,400]$}
\put(90,13){\footnotesize$I_{4} = [1,\,2,\,\dots,\,100]$, $I_{5} = [101,\,102,\,\dots,\,200]$, \dots}
\put(90, 2){\footnotesize$I_{8} = [1,\,2,\,\dots,\,50]$, $I_{9} = [51,\,52,\,\dots,\,100]$, \dots}
\put(52.5,36.5){$1$}
\put(35,25){$2$}
\put(70,25){$3$}
\put(26,13){$4$}
\put(44,13){$5$}
\put(61,13){$6$}
\put(78.5,13){$7$}
\put(22, 2){\small$8$}
\put(31, 2){\small$9$}
\put(38.5, 2){\small$10$}
\put(47, 2){\small$11$}
\put(56, 2){\small$12$}
\put(64.5, 2){\small$13$}
\put(73, 2){\small$14$}
\put(82, 2){\small$15$}
\end{picture}}
\caption{Numbering of nodes in a fully populated binary tree with $L=3$ levels.
The root is the original index vector $I = I_{1} = [1,\,2,\,\dots,\,400]$.}
\label{fig:tree}
\end{figure}

\subsection{Hierarchically Block Separable (HBS) matrices}
\label{sec:HBS}
We say that a dense $N\times N$ matrix $\mtx{A}$ is HBS with respect to a given tree
partitioning of its index vector $I = [1,\,2,\,3,\,\dots,\,N]$ if for every node $\tau$ in the tree,
there exist basis matrix $\mhU_{\tau}$ and $\mhV_{\tau}$ with the following properties:
\begin{enumerate}
\item For every sibling pair $\{\alpha,\beta\}$, the corresponding
off-diagonal block $\mtx{A}(I_{\alpha},I_{\beta})$ admits (up to precision $\varepsilon$)
the factorization
\begin{equation}
\label{eq:sibling_factorization}
\begin{array}{ccccccc}
\mtx{A}(I_{\alpha},I_{\beta}) &=& \mhU_{\alpha} & \mtA_{\alpha,\beta} & \mhV_{\beta},\\
n_{\alpha} \times n_{\beta} && n_{\alpha} \times k_{\alpha} & k_{\alpha} \times k_{\beta} & k_{\beta} \times n_{\beta}
\end{array}
\end{equation}
where $k_{\alpha} < n_{\alpha}$ and $k_{\beta} < n_{\beta}$.
\item For every sibling pair $\{\alpha,\beta\}$ with parent $\tau$, there exist
matrices $\mU_{\tau}$ and $\mV_{\tau}$ such that
\begin{align}
\label{eq:hbases_repeat}
&\begin{array}{cccccc}
\mhU_{\tau} &=& \mtwo{\mhU_{\alpha}}{\mtx{0}}{\mtx{0}}{\mhU_{\beta}} & \mU_{\tau},\\
n_{\tau}\times k_{\tau} && n_{\tau}\times (k_{\alpha} + k_{\beta}) & (k_{\alpha} + k_{\beta}) \times k_{\tau}
\end{array}\\
&\begin{array}{cccccc}
\mhV_{\tau} &=& \mtwo{\mhV_{\alpha}}{\mtx{0}}{\mtx{0}}{\mhV_{\beta}} & \mV_{\tau}.\\
n_{\tau}\times k_{\tau} && n_{\tau}\times (k_{\alpha} + k_{\beta}) & (k_{\alpha} + k_{\beta}) \times k_{\tau}
\end{array}
\end{align}
\end{enumerate}

For notational convenience, we set for every leaf node $\tau$
$$
\mU_{\tau} = \mhU_{\tau},
\qquad\mbox{and}\qquad
\mV_{\tau} = \mhV_{\tau}.
$$
Note that every basis matrix $\mU_{\tau}$ and $\mV_{\tau}$ is small.
The point is that the matrix $\mA$ is fully specified by giving
the following for every node $\tau$:
\begin{itemize}
 \item[(i)] the two small matrices $\mU_{\tau}$, and $\mV_{\tau}$,
\item[(ii)] if $\tau$ is a leaf, the diagonal block $\mtx{A}(I_{\tau},I_{\tau})$, and
\item[(iii)] if $\tau$ is a parent node with children $\{\alpha,\beta\}$, the sibling interaction
matrices $\mtA_{\alpha,\beta},\,\mtA_{\beta,\alpha}$.
\end{itemize}
\textit{Particular attention should be paid to the fact that the long basis matrices $\mhU_{\tau}$
and $\mhV_{\tau}$ are \textbf{never} formed --- these were introduced merely to facilitate
the derivation of the HBS representation.} Table \ref{fig:summary_of_factors} summarizes
the factors required, and Algorithm I describes how for any vector $\vct{\sigma}$ the matrix-vector
product $\uu = \mtx{A}\vct{\sigma}$ can be computed if the HBS factors of $\mtx{A}$ are given.


\begin{figure}
\begin{tabular}{|l|l|l|l|} \hline
                             & Name:            & Size:       & Function: \\ \hline
For each leaf                & $\mtx{A}(I_{\tau},I_{\tau})$ & $ n_{\tau}\times  n_{\tau}$ & Diagonal block. \\
node $\tau$:                 & $\mtx{U}_{\tau}$ & $ n_{\tau}\times  k_{\tau}$ & Basis for the columns in $\mA(I_{\tau},I_{\tau}^{\rm c})$. \\
                             & $\mtx{V}_{\tau}$ & $ n_{\tau}\times  k_{\tau}$ & Basis for the rows    in $\mA(I_{\tau}^{\rm c},I_{\tau})$. \\ \hline
For each parent              & $\mtA_{\alpha,\beta}$ & $k_{\alpha}\times k_{\beta}$ & Sibling interaction matrix. \\
node $\tau$ with             & $\mtA_{\beta,\alpha}$ & $k_{\beta}\times k_{\alpha}$ & Sibling interaction matrix. \\
children $\{\alpha,\beta\}$: & $\mtx{U}_{\tau}$ & $(k_{\alpha}+k_{\beta})\times  k_{\tau}$ & Basis for the (reduced) incoming fields on $\tau$. \\
                             & $\mtx{V}_{\tau}$ & $(k_{\alpha}+k_{\beta})\times  k_{\tau}$ & Basis for the (reduced) outgoing fields from $\tau$. \\ \hline
\end{tabular}
\caption{An HBS matrix $\mA$ associated with a tree $\mathcal{T}$ is fully specified if the factors
listed above are provided.}
\label{fig:summary_of_factors}
\end{figure}

\begin{figure}
\begin{center}
\fbox{
\begin{minipage}{.9\textwidth}

\begin{center}
\textsc{Algorithm I: HBS matrix-vector multiply}
\end{center}

\lsp

\textit{Given a vector $\vct{\sigma}$ and a matrix $\mA$ in HBS format, compute $\vct{u} = \mA\,\vct{\sigma}$.\\
It is assumed that the nodes are ordered so that if $\tau$ is the parent of $\sigma$, then $\tau < \sigma$.}

\lsp

\begin{tabbing}
\hspace{5mm} \= \hspace{5mm} \= \hspace{5mm} \= \kill
\textbf{for} $\tau = N_{\rm boxes},\,N_{\rm boxes}-1,\,\dots,\,2$\\
\> \textbf{if} $\tau$ is a leaf\\
\> \> $\tilde{\vct{\sigma}}_{\tau} = \mV_{\tau}^{*}\,\vct{\sigma}(I_{\tau})$.\\
\> \textbf{else}\\
\> \> $\tilde{\vct{\sigma}}_{\tau} = \mV_{\tau}^{*}\,\vtwo{\tilde{\vct{\sigma}}_{\alpha}}{\tilde{\vct{\sigma}}_{\beta}}$, where $\alpha$ and $\beta$ denote the children of $\tau$.\\
\> \textbf{end if}\\
\textbf{end for}\\
\\
$\tilde{\vct{u}}_{1} = 0$\\
\textbf{for} $\tau = 1,\,2,\,3,\,\dots,\,N_{\rm boxes}$\\
\> \textbf{if} $\tau$ is a parent\\
\> \> $\vtwo{\tilde{\vct{u}}_{\alpha}}{\tilde{\vct{u}}_{\beta}} =
       \mU_{\tau}\,\tilde{\vct{u}}_{\tau} +
       \mtwo{\mtx{0}}{\mtA_{\alpha,\beta}}{\mtA_{\beta,\alpha}}{\mtx{0}}\,
       \vtwo{\tilde{\vct{\sigma}}_{\alpha}}{\tilde{\vct{\sigma}}_{\beta}}$,
       where $\alpha$ and $\beta$ denote the children of $\tau$.\\
\> \textbf{else} \\
\> \> $\vct{u}(I_{\tau}) = \mU_{\tau}\,\tilde{\vct{u}}_{\tau} + \mA(I_{\tau},I_{\tau})\,\vct{\sigma}(I_{\tau})$.\\
\> \textbf{end if}\\
\textbf{end for}
\end{tabbing}

\end{minipage}}
\end{center}
\end{figure}

\begin{remark}[Meaning of basis matrices $\mhU_{\tau},\,\mhV_{\tau},\,\mU_{\tau},\,\mV_{\tau}$]
Let us describe the heuristic meaning of the ``tall thin'' basis matrix $\mhU_{\tau}$.
Conditions (\ref{eq:sibling_factorization}) and (\ref{eq:hbases_repeat}) together imply that
the columns of $\mhU_{\tau}$ span the column space of the off-diagonal block $\mtx{A}(I_{\tau},I_{\tau}^{\rm c})$,
as well as the columns of $\mtx{A}_{\rm in}(I_{\tau},:)$. We can write the restriction of the local
equilibrium equation on $\Gamma_{\tau}$ as
\begin{equation}
\label{eq:tau_equilibrium}
\mtx{A}(I_{\tau},I_{\tau})\,\vct{\sigma}(I_{\tau}) =
\mtx{A}_{\rm in}(I_{\tau},\colon)\,\vct{s}_{\rm in} -
\mtx{A}(I_{\tau},I_{\tau}^{\rm c})\,\vct{\sigma}(I_{\tau}^{\rm c})
\end{equation}
Exploiting that the columns of $\mhU_{\tau}$ span the column spaces of both matrices on the right hand
side of (\ref{eq:tau_equilibrium}), we can rewrite the equation as
\begin{equation}
\label{eq:tau_equilibrium2}
\mtx{A}(I_{\tau},I_{\tau})\,\vct{\sigma}(I_{\tau}) = \mhU_{\tau}\,\bigl(\tilde{\vct{u}}_{\tau} - \tilde{\vct{v}}_{\tau}\bigr),
\end{equation}
where
$\tilde{\vct{u}}_{\tau} = \mhU_{\tau}^{\dagger}\,\mtx{A}_{\rm in}(I_{\tau},\colon)\,\vct{s}_{\rm in}$
is an efficient representation of the local incoming field, and where
$\tilde{\vct{v}}_{\tau} = \mhU_{\tau}^{\dagger}\,\mtx{A}(I_{\tau},I_{\tau}^{\rm c})\,\vct{\sigma}(I_{\tau}^{\rm c})$
is an efficient representation of the  field on $\Gamma_{\tau}$ caused by charges on $\Gamma \backslash \Gamma_{\tau}$.
The basis matrix $\mhV_{\tau}$ analogously provides an efficient way to represent the outgoing fields on $\Gamma_{\tau}$.
The ``small'' matrices $\mU_{\tau}$ and $\mV_{\tau}$ then provide compact representations of $\mhU_{\tau}$ and $\mhV_{\tau}$.
Observe for instance that
$$
\mhU_{1}
=
\left[\begin{array}{cc}
\mhU_{2} & \mtx{0}  \\
\mtx{0}  & \mhU_{3}
\end{array}\right]\,
\mU_{1}
=
\left[\begin{array}{cccc}
\mhU_{4} & \mtx{0}  & \mtx{0}  & \mtx{0}  \\
\mtx{0}  & \mhU_{5} & \mtx{0}  & \mtx{0}  \\
\mtx{0}  & \mtx{0}  & \mhU_{6} & \mtx{0}  \\
\mtx{0}  & \mtx{0}  & \mtx{0}  & \mhU_{7}
\end{array}\right]\,
\left[\begin{array}{cc}
\mU_{2} & \mtx{0}  \\
\mtx{0} & \mU_{3}
\end{array}\right]\,
\mU_{1}.
$$
These concepts are described in further detail in Appendix \ref{app:fields}.
\end{remark}

\begin{remark}(Exact versus numerical rank)
In practical computations, the off-diagonal blocks of discretized integral operators
only have low \textit{numerical} (as opposed to exact) rank. However, in the present
context, the singular values tend to decay exponentially fast, which means that the
truncation error can often be rendered close to machine precision. We do not in this
manuscript provide a rigorous error analysis, but merely observe that across a broad
range of numerical experiments, we did not see any error aggregation; instead, the local
truncation errors closely matched the global error.
\end{remark}

\subsection{Hierarchical construction of a scattering matrix}
\label{sec:mergeformula}
For any node $\tau$, define its associated scattering matrix via
\begin{equation}
\label{eq:mergelemma_def1}
\mtx{S}_{\tau} = \hat{\mtx{V}}_{\tau}^{*}\,\bigl(\mtx{A}(I_{\tau},I_{\tau})\bigr)^{-1}\,\hat{\mtx{U}}_{\tau}.
\end{equation}
The following lemma states that the scattering matrix for a parent node $\tau$ with children $\alpha$
and $\beta$ can be formed inexpensively from the scattering matrices $\mtx{S}_{\alpha}$ and $\mtx{S}_{\beta}$,
along with the sibling interaction matrices $\mtA_{\alpha,\beta}$ and $\mtA_{\beta,\alpha}$.

\begin{lemma}
\label{lemma:mergelemma}
Let $\tau$ be a node with children $\alpha$ and $\beta$. Then
\begin{equation}
\label{eq:mergelemmaA}
\mtx{S}_{\tau} =
\mtx{V}_{\tau}^{*}
\mtwo{\mtx{I}}{\mtx{S}_{\alpha}\mtA_{\alpha,\beta}}{\mtx{S}_{\beta}\mtA_{\beta,\alpha}}{\mtx{I}}^{-1}
\mtwo{\mtx{S}_{\alpha}}{\mtx{0}}{\mtx{0}}{\mtx{S}_{\beta}}
\mtx{U}_{\tau}.
\end{equation}
\end{lemma}

Lemma \ref{lemma:mergelemma} immediately leads to an $O(N^{3/2})$ algorithm for computing the
scattering matrix $\mS_{1}$. A pseudo-code description is given as Algorithm II.

\vspace{1.5mm}

\noindent
\textit{\textbf{Proof of Lemma \ref{lemma:mergelemma}:}} Using formula (\ref{eq:sibling_factorization}) for the factorization of sibling interaction matrices we get
\begin{align}
\nonumber
\mtx{A}(I_{\tau},I_{\tau}) =&\
\mtwo{\mtx{A}_{\alpha,\alpha}}{\mtx{A}_{\alpha,\beta}}{\mtx{A}_{\beta,\alpha}}{\mtx{A}_{\beta,\beta}}
=
\mtwo{\mtx{A}_{\alpha,\alpha}}{\mtx{0}}{\mtx{0}}{\mtx{A}_{\beta,\beta}}
+
\mtwo{\mtx{0}}{\hat{\mtx{U}}_{\alpha}\mtA_{\alpha,\beta}\hat{\mtx{V}}_{\beta}^{*}}{\hat{\mtx{U}}_{\beta}\mtA_{\beta,\alpha}\hat{\mtx{V}}_{\alpha}^{*}}{\mtx{0}} \\
\nonumber
=&\
\mtwo{\mtx{A}_{\alpha,\alpha}}{\mtx{0}}{\mtx{0}}{\mtx{A}_{\beta,\beta}}
+
\mtwo{\mtx{0}}{\hat{\mtx{U}}_{\alpha}\mtA_{\alpha,\beta}}{\hat{\mtx{U}}_{\beta}\mtA_{\beta,\alpha}}{\mtx{0}}
\mtwo{\hat{\mtx{V}}_{\alpha}^{*}}{\mtx{0}}{\mtx{0}}{\hat{\mtx{V}}_{\beta}^{*}}\\
\nonumber
=&\
\left(
\mtx{I}
+
\mtwo{\mtx{0}}{\hat{\mtx{U}}_{\alpha}\mtA_{\alpha,\beta}}{\hat{\mtx{U}}_{\beta}\mtA_{\beta,\alpha}}{\mtx{0}}
\mtwo{\hat{\mtx{V}}_{\alpha}^{*}\mtx{A}_{\alpha,\alpha}^{-1}}{\mtx{0}}{\mtx{0}}{\hat{\mtx{V}}_{\beta}^{*}\mtx{A}_{\beta,\beta}^{-1}}
\right)
\mtwo{\mtx{A}_{\alpha,\alpha}}{\mtx{0}}{\mtx{0}}{\mtx{A}_{\beta,\beta}}\\
\label{eq:prewood}
=&\
\bigl(\mtx{I} + \mtx{E}\,\mtx{F}^{*}\bigr)\,
\mtwo{\mtx{A}_{\alpha,\alpha}}{\mtx{0}}{\mtx{0}}{\mtx{A}_{\beta,\beta}},
\end{align}
where we defined
$$
\mtx{E} =
\mtwo{\mtx{0}}{\hat{\mtx{U}}_{\alpha}\mtA_{\alpha,\beta}}{\hat{\mtx{U}}_{\beta}\mtA_{\beta,\alpha}}{\mtx{0}}
\qquad\mbox{and}\qquad
\mtx{F}^{*} =
\mtwo{\hat{\mtx{V}}_{\alpha}^{*}\mtx{A}_{\alpha,\alpha}^{-1}}{\mtx{0}}{\mtx{0}}{\hat{\mtx{V}}_{\beta}^{*}\mtx{A}_{\beta,\beta}^{-1}}.
$$
Formula (\ref{eq:prewood}) shows that $\mtx{A}(I_{\tau},I_{\tau})$ can be written as a product between a low-rank perturbation to the identity,
and a block-diagonal matrix. The Woodbury formula for inversion of a low-rank perturbation of the identity states that
$(\mtx{I} + \mtx{E}\mtx{F}^{*})^{-1} = \mtx{I} - \mtx{E}(\mtx{I}+\mtx{F}^{*}\mtx{E})^{-1}\mtx{F}^{*}$.
To exploit this formula, we first define
\begin{equation}
\label{eq:autistic11x}
\mtx{Z}_{\tau} \stackrel{\mbox{def}}{=}
\big(\mtx{I} + \mtx{F}^{*}\mtx{E}\bigr)^{-1} =
\left(\mtx{I} + \mtwo{\mtx{0}}{\hat{\mtx{V}}_{\alpha}^{*}\mtx{A}_{\alpha,\alpha}^{-1}\hat{\mtx{U}}_{\alpha}\,\mtA_{\alpha,\beta }}
                              {\hat{\mtx{V}}_{\beta }^{*}\mtx{A}_{\beta ,\beta }^{-1}\hat{\mtx{U}}_{\beta }\,\mtA_{\beta ,\alpha}}{0}\right)^{-1} =
\mtwo{\mtx{I}}{\mtx{S}_{\alpha}\mtA_{\alpha,\beta}}{\mtx{S}_{\beta}\mtA_{\beta,\alpha}}{\mtx{I}}^{-1}.
\end{equation}
Then
\begin{equation}
\label{eq:autistic12x}
\mtx{A}(I_{\tau},I_{\tau})^{-1} =
\mtwo{\mtx{A}_{\alpha,\alpha}^{-1}}{\mtx{0}}{\mtx{0}}{\mtx{A}_{\beta,\beta}^{-1}}
\left(\mtx{I} -
\mtwo{\mtx{0}}{\hat{\mtx{U}}_{\alpha}\mtA_{\alpha,\beta}}{\hat{\mtx{U}}_{\beta}\mtA_{\beta,\alpha}}{\mtx{0}}
\mtx{Z}_{\tau}
\mtwo{\hat{\mtx{V}}_{\alpha}^{*}\mtx{A}_{\alpha,\alpha}^{-1}}{\mtx{0}}{\mtx{0}}{\hat{\mtx{V}}_{\beta}^{*}\mtx{A}_{\beta,\beta}^{-1}}
\right).
\end{equation}
Combining (\ref{eq:mergelemma_def1}) with (\ref{eq:autistic12x}) and inserting the condition on nestedness of the basis
matrices (\ref{eq:hbases_repeat}), we get
\begin{align*}
\mtx{S}_{\tau}
=&\
\mtx{V}_{\tau}^{*}
\mtwo{\hat{\mtx{V}}_{\alpha}^{*}}{\mtx{0}}{\mtx{0}}{\hat{\mtx{V}}_{\beta}^{*}}
\mtwo{\mtx{A}_{\alpha,\alpha}^{-1}}{\mtx{0}}{\mtx{0}}{\mtx{A}_{\beta,\beta}^{-1}}
\left(\mtx{I} -
\mtwo{\mtx{0}}{\hat{\mtx{U}}_{\alpha}\mtA_{\alpha,\beta}}{\hat{\mtx{U}}_{\beta}\mtA_{\beta,\alpha}}{\mtx{0}}
\mtx{Z}_{\tau}
\mtwo{\hat{\mtx{V}}_{\alpha}^{*}\mtx{A}_{\alpha,\alpha}^{-1}}{\mtx{0}}{\mtx{0}}{\hat{\mtx{V}}_{\beta}^{*}\mtx{A}_{\beta,\beta}^{-1}}
\right)
\mtwo{\hat{\mtx{U}}_{\alpha}}{\mtx{0}}{\mtx{0}}{\hat{\mtx{U}}_{\beta}}
\mtx{U}_{\tau}\\
=&\
\mtx{V}_{\tau}^{*}
\left(\mtx{I} -
\mtwo{\mtx{0}}{\mtx{S}_{\alpha}\mtA_{\alpha,\beta}}{\mtx{S}_{\beta}\mtA_{\beta,\alpha}}{\mtx{0}}
\mtx{Z}_{\tau}
\right)
\mtwo{\mtx{S}_{\alpha}}{\mtx{0}}{\mtx{0}}{\mtx{S}_{\beta}}
\mtx{U}_{\tau}\\
=&\
\mtx{V}_{\tau}^{*}
\left(\mtx{Z}_{\tau}^{-1} -
\mtwo{\mtx{0}}{\mtx{S}_{\alpha}\mtA_{\alpha,\beta}}{\mtx{S}_{\beta}\mtA_{\beta,\alpha}}{\mtx{0}}
\right)
\mtx{Z}_{\tau}
\mtwo{\mtx{S}_{\alpha}}{\mtx{0}}{\mtx{0}}{\mtx{S}_{\beta}}
\mtx{U}_{\tau}
=
\mtx{V}_{\tau}^{*}
\mtx{Z}_{\tau}
\mtwo{\mtx{S}_{\alpha}}{\mtx{0}}{\mtx{0}}{\mtx{S}_{\beta}}
\mtx{U}_{\tau}.
\end{align*}
So (\ref{eq:mergelemmaA}) holds.\qed

\begin{figure}
\fbox{\begin{minipage}{\textwidth}
\begin{center}
\textsc{Algorithm II: Hierarchical computation of scattering matrices}
\end{center}

\lsp

\begin{tabbing}
\hspace{5mm} \= \hspace{5mm} \= \hspace{5mm} \= \hspace{5mm} \= \hspace{5mm} \kill
\textbf{for} $\tau = N_{\rm boxes} : (-1) : 1$\\
\> \textbf{if} $\tau$ is a leaf\\
\> \> $\mtx{S}_{\tau} = \mtx{V}_{\tau}^{*}\,\mtx{A}(I_{\tau},I_{\tau})^{-1}\,\mtx{U}_{\tau}$.\\
\> \textbf{else}\\
\> \> $\mtx{S}_{\tau} = \mtx{V}_{\tau}^{*}
                        \left[\begin{array}{cc} \mtx{I} & \mtx{S}_{\alpha}\mtA_{\alpha,\beta} \\ \mtx{S}_{\beta}\mtA_{\beta,\alpha} & \mtx{I} \end{array}\right]^{-1}
                        \left[\begin{array}{cc}\mtx{S}_{\alpha} & \mtx{0} \\ \mtx{0} & \mtx{S}_{\beta} \end{array}\right]
                        \mtx{U}_{\tau}$,
                        where $\alpha$ and $\beta$ denote the children of $\tau$.\\
\> \textbf{end if}\\
\textbf{end for}\\
\end{tabbing}
\end{minipage}}
\end{figure}

\begin{remark}[Extension to a full direct solver]
The scattering matrix formulation described in this section can be viewed as
an efficient direct solver for the equation $\mtx{A}\vct{\sigma} = \vct{u}$
for the particular case where the load vector $\vct{u}$ belongs to the low-dimensional
subspace spanned by the columns of the given matrix $\mhU_{1}$.
A general direct solver for applying $\mtx{A}^{-1}$ to an arbitrary vector
for an HBS matrix $\mtx{A}$ (under some general conditions on invertibility) can
be obtained by slight variations to the simplified scheme described here, see
Appendix \ref{app:fields} for details.
\end{remark}

\subsection{Efficient construction of the HBS representation and recursive skeletonization}
\label{sec:skeletons}

In this section, we describe a technique for computing the HBS representation
of a discretized boundary integral operator via a procedure described in
\cite{2005_martinsson_fastdirect,2009_martinsson_ACTA}, and sometimes referred to as
\textit{recursive skeletonization} \cite{2012_ho_greengard_fastdirect}.

\subsubsection{The interpolative decomposition (ID)}
The \textit{interpolative decomposition} (ID) \cite{2005_martinsson_skel} is a low-rank factorization
in which a matrix of rank $k$ is expressed by using a selection of $k$ of its columns (rows) as a basis
for its column (row) space. To be precise, let $\mtx{X}$ be matrix of size $m\times n$ and of rank $k$. Then
it is possible to determine an index vector $J_{\rm colskel} \subseteq \{1,2,3,\dots,n\}$ of length $k$
such that
\begin{equation}
\label{eq:colskel}
\begin{array}{cccccc}
\mtx{X} &=& \mtx{X}(\colon,J_{\rm colskel})&\mV^{*},\\
m\times n && m\times k & k\times n
\end{array}
\end{equation}
where $\mV$ is an $n\times k$ matrix that contains the $k\times k$ unit matrix as a sub-matrix,
$$
\mV(J_{\rm colskel},\colon) = \mtx{I}.
$$
Moreover, no entry of $\mV$ is of magnitude greater than $1$, which implies that $\mV$ is well-conditioned.
We call the index vector $J_{\rm colskel}$ the \textit{column skeleton} of $\mtx{X}$.

Computing the optimal ID of a general matrix is hard \cite{gu1996},
but if we slightly relax the restriction on the basis matrix $\mV$
to allow its entries to be bounded by, say, $2$ instead of one, then
very efficient algorithms are available~\cite{2005_martinsson_skel,gu1996,2011_martinsson_randomsurvey}.

\subsubsection{Interpolative decompositions and HBS matrices}
\label{sec:HBS_ID}
It is highly advantageous to use the ID to factor the off-diagonal blocks in
the HBS representation of a matrix.
The first step in the construction is to form an extended system matrix
$$
\mtx{B} =
\left[\begin{array}{cc}
\mtx{A} & \mtx{A}_{\rm in} \\
\mtx{A}_{\rm out} & \mtx{0}
\end{array}\right].
$$
For each leaf node $\tau$, identify an ID basis matrix $\mU_{\tau}$ and a subset
$\tilde{I}^{\rm in}_{\tau} \subset I_{\tau}$ such that
\begin{equation}
\label{eq:ram1}
\begin{array}{ccccc}
\mtx{B}(I_{\tau},I_{\tau}^{\rm c}) &=& \mU_{\tau} & \mtx{B}(\tilde{I}^{\rm in}_{\tau},I_{\tau}^{\rm c}).\\
n_{\tau} \times (N + N_{\rm in}-n_{\tau}) && n_{\tau} \times k_{\tau} & k_{\tau} \times (N + N_{\rm in}-n_{\tau})
\end{array}
\end{equation}
The vector $\tilde{I}^{\rm in}_{\tau}$ marks a subset of the collocation points in $\Gamma$ that we
call the \textit{skeleton points of $\tau$}. These points have the property that all interactions
between $\Gamma_{\tau}$ and the rest of the geometry can be done via evaluation through these points.
Next let $\tau$ be a node whose children $\{\alpha,\beta\}$ are leaves. Then we find the
skeleton vector $\tilde{I}^{\rm in}_{\tau}$ and the basis matrix $\mU_{\tau}$ by simply factoring
the matrix $\mtx{B}([\tilde{I}^{\rm in}_{\alpha},\tilde{I}^{\rm in}_{\beta}],I_{\tau}^{\rm c})$ to form its ID
\begin{equation}
\label{eq:ram2}
\begin{array}{ccccc}
\mtx{B}([\tilde{I}_{\alpha}^{\rm in},\tilde{I}_{\beta}^{\rm in}],I_{\tau}^{\rm c}) &=& \mU_{\tau} & \mtx{B}(\tilde{I}^{\rm in}_{\tau},I_{\tau}^{\rm c}).\\
(k_{\alpha}+k_{\beta}) \times (N + N_{\rm in} - n_{\tau}) && (k_{\alpha}+k_{\beta}) \times k_{\tau} & k_{\tau} \times (N + N_{\rm in} - n_{\tau})
\end{array}
\end{equation}
All remaining skeleton sets and basis matrices $\{\mU_{\tau}\}_{\tau}$ are found by simply
continuing this process up through the tree.

The process for determining the basis matrices $\{\mV_{\tau}\}_{\tau}$ is analogous. For each
leaf form the factorization
\begin{equation}
\label{eq:ram3}
\begin{array}{ccccc}
\mtx{B}(I_{\tau}^{\rm c},I_{\tau}) &=& \mtx{B}(\tilde{I}^{\rm out}_{\tau},I_{\tau}^{\rm c}) & \mV_{\tau}^{*}.\\
(N + N_{\rm out}-n_{\tau}) \times n_{\tau} && (N + N_{\rm out}-n_{\tau}) \times k_{\tau} & k_{\tau} \times n_{\tau}
\end{array}
\end{equation}
For a parent $\tau$ with children $\{\alpha,\beta\}$ (whose skeletons are available), factor
\begin{equation}
\label{eq:ram4}
\begin{array}{ccccc}
\mtx{B}(I_{\tau}^{\rm c},[\tilde{I}^{\rm out}_{\alpha},\tilde{I}^{\rm out}_{\beta}]) &=& \mtx{B}(I_{\tau}^{\rm c},\tilde{I}^{\rm out}_{\tau}) & \mV_{\tau}^{*}.\\
(N + N_{\rm out}-n_{\tau}) \times (k_{\alpha}+k_{\beta}) && (N + N_{\rm out}-n_{\tau}) \times k_{\tau} & k_{\tau} \times (k_{\alpha}+k_{\beta})
\end{array}
\end{equation}

Once the skeletons $\{\tilde{I}^{\rm in}_{\tau},\,\tilde{I}^{\rm out}_{\tau}\}$ have been determined for all nodes, the sibling interaction matrices
for a sibling pair $\{\alpha,\beta\}$ are given simply by
\begin{equation}
\label{eq:ram5}
\mtA_{\alpha,\beta} = \mtx{A}(\tilde{I}_{\alpha}^{\rm out},\tilde{I}_{\beta}^{\rm in}),
\qquad\mbox{and}\qquad
\mtA_{\beta,\alpha} = \mtx{A}(\tilde{I}_{\beta}^{\rm out},\tilde{I}_{\alpha}^{\rm in}).
\end{equation}
The formulas (\ref{eq:ram5}) are very useful since they make it extremely cheap to form the sibling interaction
matrices --- they are merely submatrices of the original system matrix!

\begin{remark}
It is often convenient to force the incoming and the outgoing skeletons to be identical,
$\tilde{I}^{\rm out}_{\tau} = \tilde{I}^{\rm in}_{\tau}$. For instance, in order to determine
a single set of skeleton nodes $\tilde{I}_{\tau}$ for a leaf $\tau$, we can combine
(\ref{eq:ram1}) and (\ref{eq:ram2}) to the single factorization
\begin{equation}
\label{eq:ram6}
\begin{array}{ccccc}
\vtwo{\mtx{B}(I_{\tau}^{\rm c},I_{\tau})}{\mtx{B}(I_{\tau},I_{\tau}^{\rm c})^{*}}
&=& \vtwo{\mtx{B}(I_{\tau}^{\rm c},\tilde{I}_{\tau})}{\mtx{B}(\tilde{I}_{\tau},I_{\tau}^{\rm c})^{*}} & \mV_{\tau}^{*}.\\
(2N + N_{\rm out} + N_{\rm in } - n_{\tau}) \times n_{\tau} && (2N + N_{\rm out} + N_{\rm in} - n_{\tau}) \times k_{\tau} & k_{\tau} \times n_{\tau}
\end{array}
\end{equation}
Then simply set $\mU_{\tau} = \mV_{\tau}$.
\end{remark}

\begin{figure}[h!]
\centering
 \begin{tabular}{ccc}
\includegraphics[width=70mm]{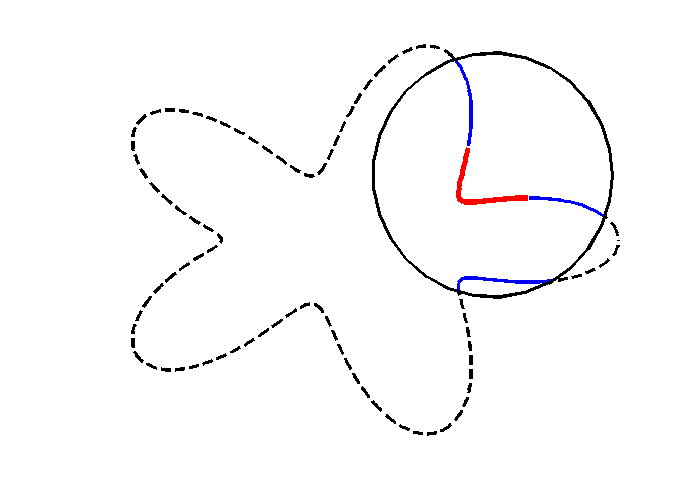}
 &
 \mbox{}
 &
\includegraphics[width=80mm]{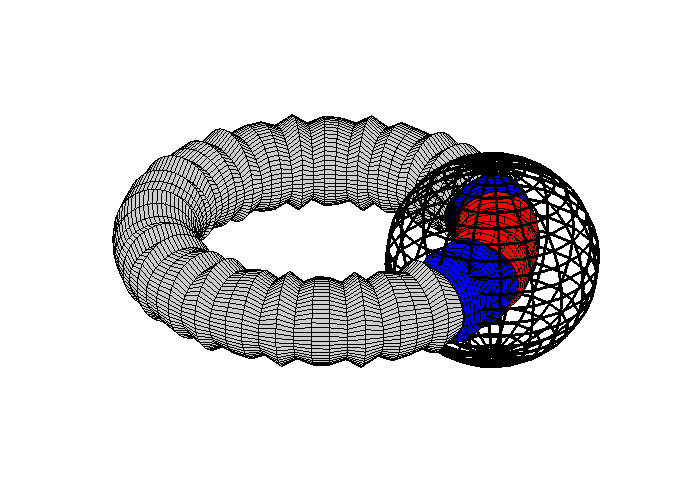}
 \\
 (a) & & (b)
 \end{tabular}
\caption{\label{fig:proxy} Illustration of the geometric objects $\Gamma_{\tau}$, $\Gamma_{\rm near}$ and
$\Gamma_{\rm proxy}$ used for the Green's representation technique for a two-dimensional (a) and
three-dimensional surface (b). $\Gamma_{\tau}$ is given in red, $\Gamma_{\rm near}$ is given in blue, and
$\Gamma_{\rm proxy}$ is given in black. }
\end{figure}

\subsubsection{The Green's representation technique}
The compression technique described in Section \ref{sec:HBS_ID} would be prohibitively
expensive if executed as stated, since it requires the formation and factorization of
very large matrices. For instance, equation (\ref{eq:ram1}) requires the factorization
of a matrix of size $n_{\tau} \times (N + N_{\rm in} 1- n_{\tau})$. It turns out that
these computations can be localized by exploiting that the column space of the matrix
$\mtx{B}(I_{\tau},I_{\tau}^{\rm c})$ (which is what we need to span) expresses merely
a set of solutions to the homogeneous Helmholtz equation. This means that we can use
representation techniques from potential theory to replace all ``far-field'' interactions
by interactions with a small proxy surface $\Gamma_{\rm proxy}$ that encloses the patch $\Gamma_{\tau}$ that
we seek to compress. All discretization points inside on $\Gamma_{\rm proxy}$ but not in $\Gamma_{\tau}$ are
labeled $\Gamma_{\rm near}$.  Figure \ref{fig:proxy} illustrates the geometry used in the Green's
representation technique for a two-dimensional contour and a three surface that is defined rigorously in
Section \ref{sec:tori}.  This technique is described in detail in Section 6.2 of
\cite{2012_martinsson_FDS_survey}. For the experiments in this paper, it is sufficient
to let $\Gamma_{\rm proxy}$ be a sphere with a radius that is twice the radius of the smallest ball containing
$\Gamma_{\tau}$. (In other words, $\Gamma_{\rm proxy}$ is the surface of the ball ``$2B$'' defined in
Section \ref{sec:quadrature}.)

\end{section}

\begin{section}{Numerical experiments.}
\label{sec:numerics}
We now present the results of a number of experiments conducted to measure the performance
of the approach of this paper.
All code was written in Fortran~77 and compiled with the
Intel Fortran Compiler version 12.1.  The experiments were carried
out on a workstation equipped with $12$ Intel Xeon processor cores
running at $3.47$ GHz and $192$ GB of RAM.   As a point of reference, we note that MATLAB
requires approximately 35 seconds to invert a $8192 \times 8192$ matrix
with complex-valued i.i.d. Gaussian entries on this workstation.


%

\begin{subsection}{Sound-soft scattering from grids of ellipsoids}

The purpose of this first set of experiments was to measure the
growth in the running time of the algorithm of Section~\ref{sec:scatteringmatrices}
as the geometry of the scatterer becomes more complicated and more nodes are
required to discretize it.  To that end we considered the exterior
Dirichlet problem (\ref{formulations:dirichlet}) on
a collection of domains consisting of grids of ellipsoids of various sizes
and eccentricities.
In each experiment, the ellipsoids were arranged on an $n\times n\times n$
grid with center points
\begin{equation*}
(1.25 + 2.5 i, 1.25 + 2.5 j, 1.25 + 2.5 k),
\qquad 1 \leq i,j,k \leq n.
\end{equation*}

\begin{figure}[!ht]
\includegraphics[width=.45\textwidth]{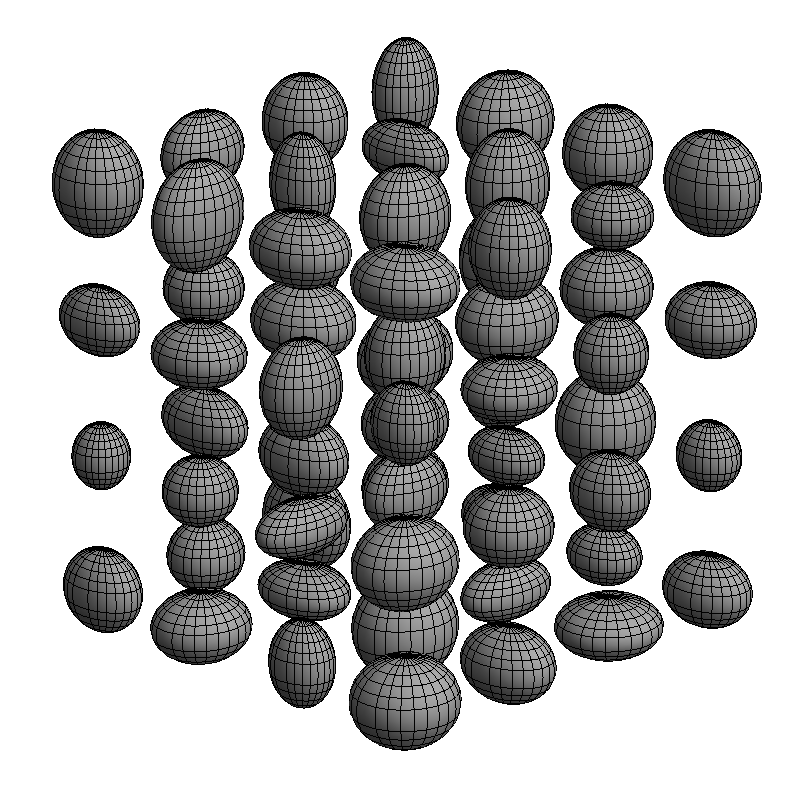}
\hfil
\includegraphics[width=.45\textwidth]{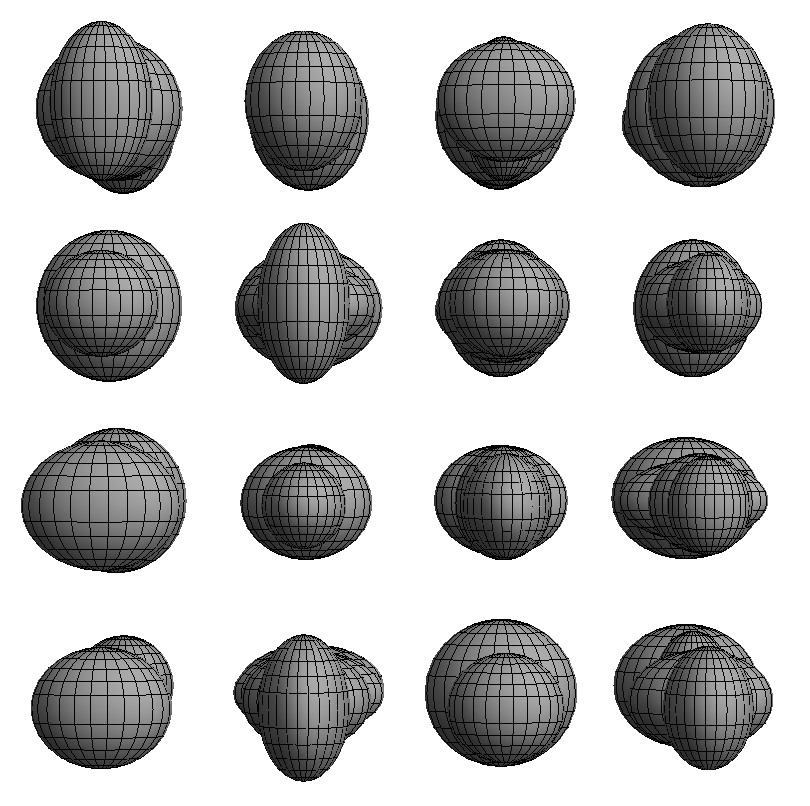}

\caption{Two views of one of the ``grid of ellipsoids'' domains of Section~\ref{section:ellipsoid_party}.  Domains
of this type pose significant challenges for iterative solvers.}
\label{figure:ellipsoid_party}
\end{figure}

\begin{table}[h!!!!!!!!!!]
\small
\begin{tabular}{crccccrr}
Grid dimensions & $N$ & $N_{\mbox{\tiny out}} \times N_{\mbox{\tiny in}}$ & $T$ & $E$ & Ratio & Predicted\\
\midrule
\addlinespace[.25em]
$2\times 2 \times 2$    & $12\,288$   &   $410\times 401$ & $1.02\e{+1}$  & $3.37\e{-04}$ & -    & -    \\
$3\times 3 \times 3$    & $41\,472$   &   $464\times 475$ & $3.43\e{+1}$  & $4.81\e{-04}$ & 3.4  & 6.2  \\
$4\times 4 \times 4$    & $98\,304$   &   $483\times 532$ & $7.92\e{+1}$  & $1.57\e{-04}$ & 2.3  & 3.7  \\
$6\times 6 \times 6$    & $331\,776$  &   $504\times 707$ & $2.96\e{+2}$  & $7.03\e{-04}$ & 3.7  & 6.2  \\
$8\times 8 \times 8$    & $786\,432$  &   $513\times 1014$ & $6.70\e{+2}$ & $4.70\e{-04}$ & 2.3  & 3.7 \\
$10\times 10 \times 10$ & $1\,536\,000$ &   $518\times 1502$ & $2.46\e{+3}$ & $3.53\e{-04}$ & 3.7  & 2.7 \\
\addlinespace[.25em]
\end{tabular}
\vskip .5em
\caption{Results for the experiments of Section~\ref{section:ellipsoid_party} in the case $\epsilon = 1.0\e{-3}$.}
\label{table:ellipsoid_party1}
\vskip -1em
\end{table}
The axis lengths of the ellipsoids were determined randomly --- each axis length was
taken to be $0.5 + 0.5 \eta$ where $\eta$ is a computer-generated pseudorandom number
between $0$ and $1$.  A $4 \times 4 \times 4$
domain of this type is depicted in Figure~\ref{figure:ellipsoid_party}.
For each experiment, the wavenumber $\kappa$ was chosen to be
\begin{equation*}
\kappa = \frac{8\pi}{2.5 \times n }
\end{equation*}
so that each domain considered was bounded by a box whose sides were
4 wavelengths in size.

Each experiment consisted of forming a scattering matrix for the
integral equation (\ref{formulations:combined}) and using it to solve an
instance of the boundary value problem (\ref{formulations:dirichlet}).
The boundary data for each problem was taken to be the restriction
of the function
\begin{equation*}
g(x) = \sum_{j=1}^N \frac{\exp(i\kappa|x-x_j|)}{|x-x_j|},
\end{equation*}
where $N$ is the total number of ellipsoids and the $x_1,\ldots,x_N$ are the
centers of the ellipsoids, to the boundary of the ellipsoids.
Let $r$ and $x_0$ denote the radius and center of the minimum bounding sphere containing
the boundary of the domain, respectively.
Then the surface $\Gamma_{\mbox{\tiny in}}$ which specified incoming potentials
was the union of the spheres of radius $0.1$ centered at the points
$x_1,\ldots,x_N$ and the sphere of radius $2r$ centered at $x_0$ while the
surface $\Gamma_{\mbox{\tiny out}}$ specifying outgoing potentials was taken to be
the sphere of radius $2r$ centered at $x_0$.
The function $g(x)$ is, of course, the unique solution of the boundary value problem
and we measured the accuracy of each obtained solution $\sigma$ of the the integral
equation by computing the relative error
\begin{equation}
E = \left(\dblint_\Gamma \left| g(x)\right|^2 dx\right)^{-1/2}
\left(\dblint_\Gamma \left| g(x) - (S_\kappa + i\kappa D_\kappa) \sigma(x)\right|^2 dx\right)^{1/2},
\label{ellipsoid_party:error}
\end{equation}
where $\Gamma$ is the sphere of radius $2$ centered at the point $x_0-(0,0,-2r-10)$.
The ellipsoids were parameterized by projection onto the cube
$[-1,1]^3$ and the faces of the cube were triangulated in order to produce
a discretization.
Three sets of experiments were performed, each at a different level of desired
precision and the ellipsoids were discretized as follows in each case:

\begin{enumerate}[$\cdot$ ]

\item
For precision $\epsilon = 1.0\e{-3}$, each cube face comprising the
parameterization domain was divided into
$8$ triangles and a $32$ point $6$th order quadrature was applied to each triangle.
This yields a total of $1536$ discretization nodes per ellipsoid.
\vskip 1em

\item
For precision $\epsilon = 1.0\e{-6}$, each cube face was subdivided
into $32$ triangles and a $32$ point $6$th order quadrature was applied to
each triangle.  This resulted in $6144$ discretization nodes per ellipsoid.
\vskip 1em

\item
For precision $\epsilon = 1.0\e{-9}$, each cube face was subdivided
into $8$ triangles and a $112$ point $12$th order quadrature was applied
to each triangle.  This results in $5376$ discretization nodes per ellipsoid.
\vskip 1em
\end{enumerate}

\begin{table}[h!!!!!!!!!!]
\small
\begin{tabular}{crccccrr}
Grid dimensions & $N$ & $N_{\mbox{\tiny out}} \times N_{\mbox{\tiny in}}$ & $T$ & $E$ & Ratio & Predicted\\
\midrule
\addlinespace[.25em]
$2\times 2 \times 2$    & $49\,152$    &  $601\times 584$ & $1.61\e{+2}$ & $1.22\e{-07}$      & -      &  -  \\
$3\times 3 \times 3$    & $165\,888$   &  $676\times 677$ & $6.87\e{+2}$ & $4.92\e{-07}$      & 4.3    & 6.2 \\
$4\times 4 \times 4$    & $393\,216$   &  $703\times 747$ & $1.68\e{+3}$ & $5.31\e{-07}$      & 2.4    & 3.6 \\
$6\times 6 \times 6$    & $1\,327\,104$  &  $728\times 925$ & $6.66\e{+3}$ & $4.60\e{-06}$      & 4.0    & 6.2 \\
$8\times 8 \times 8$    & $3\,145\,728$  & $742\times 1\,237$ & $1.59\e{+4}$ & $2.30\e{-07}$      & 2.4    & 3.6 \\
\addlinespace[.25em]
\end{tabular}
\vskip .5em
\caption{Results for the numerical experiments of Section~\ref{section:ellipsoid_party} in the case $\epsilon = 1.0\e{-6}.$ }
\label{table:ellipsoid_party2}
\vskip -1em
\end{table}
%
Tables \ref{table:ellipsoid_party1}, \ref{table:ellipsoid_party2}, and \ref{table:ellipsoid_party3}
present the results of the experiments with $\epsilon = 1.0\e{-3},\,1.0\e{-6},\,1.0\e{-9}$, respectively.
In each of those tables, $N$ refers to the
number of discretization nodes on the surface; $N_{\mbox{\tiny out}} \times N_{\mbox{\tiny in}}$
denotes the dimension of the resulting scattering matrix;
$T$ is the time
in seconds required to construct the scattering matrix;
$E$ is the error (\ref{ellipsoid_party:error});
Ratio refers to the ratio of running times; and Predicted is the
ratio of running times one would expect for an algorithm whose cost
scales as $N^{1.5}$.

\begin{table}[h!!!!!!!!!!]
\small
\begin{tabular}{crccccrr}
Grid dimensions & $N$ & $N_{\mbox{\tiny out}} \times N_{\mbox{\tiny in}}$ & $T$ & $E$ & Ratio & Predicted\\
\midrule
\addlinespace[.25em]
$2\times 2 \times 2$    & $43\,008$    &  $997\times 946$   & $5.14\e{+2}$ & $8.98\e{-09}$      & -      & -   \\
$3\times 3 \times 3$    & $145\,152$   &  $1\,363\times 1\,301$ & $2.32\e{+3}$ & $4.97\e{-10}$      & 4.5    & 6.2 \\
$4\times 4 \times 4$    & $344\,064$   &  $1\,377\times 1\,353$ & $5.56\e{+3}$ & $1.10\e{-10}$      & 2.4    & 3.6 \\
\addlinespace[.25em]
\end{tabular}
\vskip .5em
\caption{Results for the numerical experiments of Section~\ref{section:ellipsoid_party} in the case $\epsilon = 1.0\e{-9}.$ }
\label{table:ellipsoid_party3}
\vskip -2em
\end{table}

\begin{remark}
\label{re:gmres}
The boundary value problems solved in the experiments of this section present
considerable difficulties for iterative solvers.  We attempted to compare
the running time of the fast solver of this paper with the running
time of an iterative solver.  However, it took more than $1000$
iterations for GMRES residuals to fall below $1.0\e{-3}$ when we attempted
to solve the boundary value problem (\ref{formulations:dirichlet}) given on
the $2 \times 2 \times 2$ grid of ellipsoids with $1536$ discretization
nodes per ellipsoid via the combined field integral equation.
Even with the aid of an extremely efficient
multipole code, replicating the experiments of this section with an
iterative solver appears to be prohibitively expensive.
\end{remark}

\label{section:ellipsoid_party}
\end{subsection}

%
%

\begin{subsection}{Sound-hard scattering from a deformed torus.}
\label{sec:tori}
In the experiments described here we considered the problem (\ref{formulations:neumann})
with $\Gamma$ the surface parameterized via
\begin{equation}
\begin{split}
r(s,t) =
\left(
\begin{array}{c}
(2 + 0.50 \cos(s) ) \sin(t) \\
(2 + 0.50 \cos(s) ) \cos(t) \\
0.50 \sin(s) \left(1 + 0.15\cos( 36t)\right)
\end{array}
\right),
\end{split}
\ \ \ \ \ \
0 \leq s,t < 2\pi.
\ \ \
\end{equation}
Three views of this surface,
which is enclosed in the box with lower left corner $(-2.5,-2.5,-0.65)$
and upper right corner $(2.5,2.5,0.65)$, are shown in Figure~\ref{figure:deformed_torus}.

\begin{figure}[h!!]
\includegraphics[width=.32\textwidth]{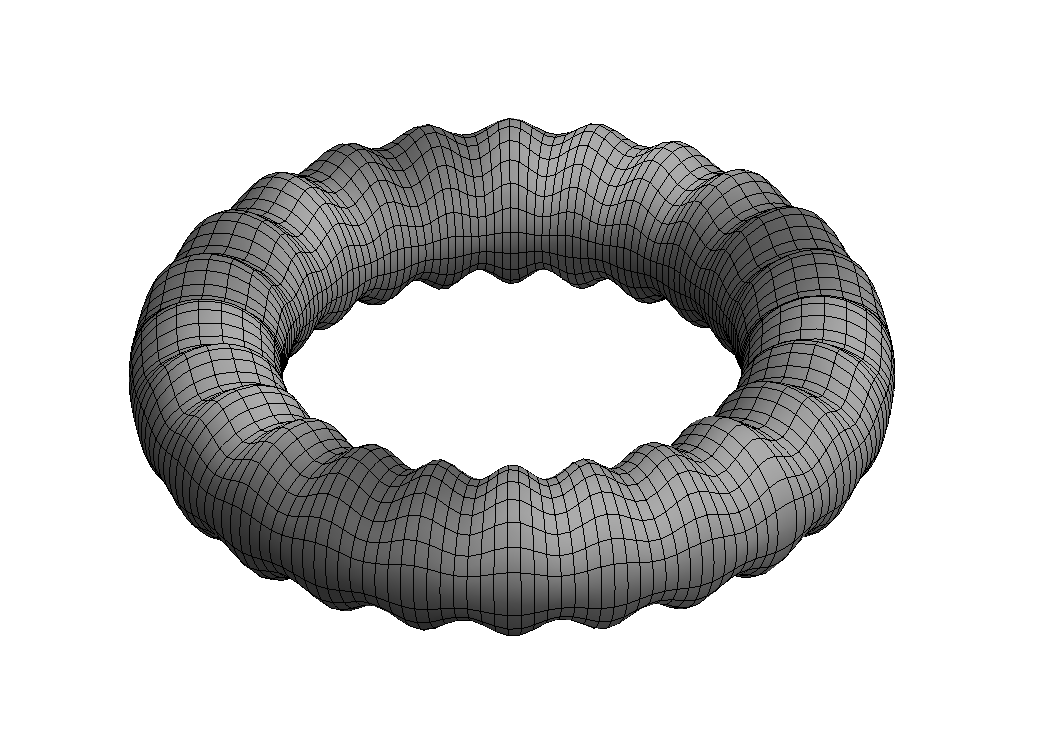}
\hfil
\raisebox{.5in}{\includegraphics[width=.32\textwidth]{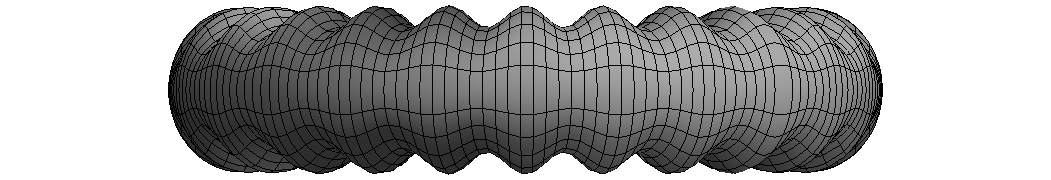}}
\hfil
\includegraphics[width=.32\textwidth]{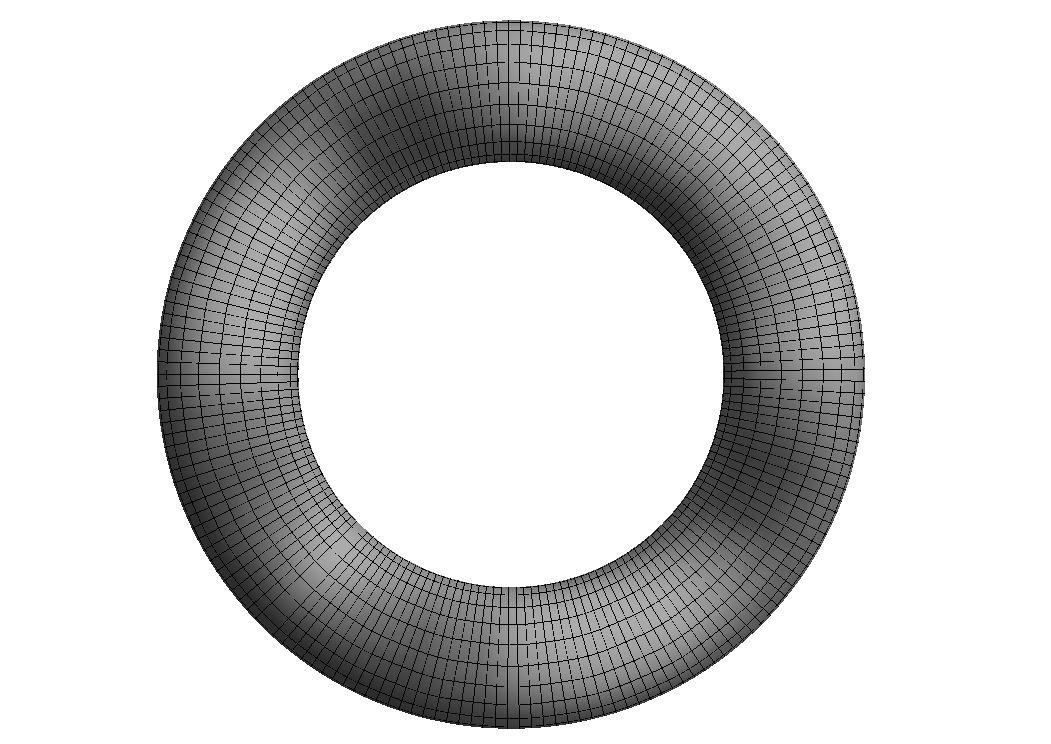}
\caption{Three views of the deformed torus of Section~\ref{section:deformed_torus}.}
\label{figure:deformed_torus}
\end{figure}

Each experiment consisted of constructing a scattering matrix
 for the integral equation formulation (\ref{formulations:naive})
of the problem (\ref{formulations:neumann}) and using that scattering matrix
to solve an instance of the boundary value problem.
The wavenumber $\kappa$ was set at $4\pi/5$, which makes the domain
approximately $2 \times 2 \times 0.7$ wavelengths in size,
and the boundary data was taken to be the normal derivative of the function
\begin{equation*}
g(x) = \frac{1}{4\pi} \frac{\exp\left(i \kappa \left|x-x_0\right|\right)}{\left|x-x_0\right|},
\end{equation*}
where $x_0$ is the point $(2,0,0)$.  Since
$x_0$ lies inside of the domain $\Omega$, the function $g(x)$ is, of course, the
unique solution of (\ref{formulations:neumann}).
The surface $\Gamma_{\mbox{\tiny in}}$ which specified incoming potentials
was take to be the union of the sphere of radius $0.1$ centered at the point
$(2,0,0)$ and the sphere of radius $6$ centered at $(0,0,0)$ while the
surface $\Gamma_{\mbox{\tiny out}}$ specifying outgoing potentials was taken to be
the sphere of radius $6$ centered at $(0,0,0)$.
The accuracy of each obtained solution $\sigma$ was measured by computing
the relative error
\begin{equation}
E =
\left(\dblint_S \left| g(x)\right|^2 dx\right)^{-1/2}
\left(\dblint_S \left| g(x) - S_\kappa \sigma(x)\right|^2 dx\right)^{1/2},
\label{deformed_torus:error}
\end{equation}
where $S$ is the sphere of radius $1$ centered at the point $(-20,0,0)$.
The order of the discretization quadrature was fixed at $8$ and the requested
precision for the scattering matrix was set to $\epsilon = 1.0\e{-7}$ in
each experiment.  The number of triangles used to decompose the parameterization
domain was varied.

\begin{table}[h!!]
\small
\begin{tabular}{rrccccc}
\addlinespace[.5em]
$N_{\mbox{\tiny tris}}$ & $N$ & $T$ & $E$  & $N_{\mbox{\tiny out}} \times N_{\mbox{\tiny in}}$ & Ratio &
Predicted\\
\midrule
\addlinespace[.5em]
$32$   &   $1\,664$    & $7.16\e{+00}$ & $3.51\e{-02}$ & $395\times 447$ & -      & -\\
$128$  &   $6\,656$    & $6.29\e{+01}$ & $4.41\e{-03}$ & $407\times 460$ & $8.79$ & $8$ \\
$512$  &   $26\,624$   & $2.81\e{+02}$ & $4.08\e{-05}$ & $400\times 456$ & $4.47$ & $8$\\
$2\,048$ &   $106\,496$  & $2.60\e{+03}$ & $7.80\e{-07}$ & $408\times 463$ & $9.25$ & $8$\\
$8\,192$ &   $425\,984$  & $1.47\e{+04}$ & $3.25\e{-08}$ & $409\times 466$ & $5.66$ & $8$\\
\end{tabular}
\vskip 1em
\caption{The results of the experiments of Section~\ref{section:deformed_torus}.}
\label{table:deformed_torus}
\vskip -2em
\end{table}

Table~\ref{table:deformed_torus} shows the results.
There, $N_{\mbox{\tiny tris}}$ is the number of triangles
into which the parameterization domain was partitioned;
 $N$ is the total number of discretization nodes on the boundary $\Gamma$;
$T$ is the total running time of the algorithm in seconds;
$E$ is the relative error (\ref{deformed_torus:error});
$N_{\mbox{\tiny out}}\times N_{\mbox{\tiny in}}$ gives the dimensions of the
scattering matrix; Ratio is the quotient of the running time at the current
level of discretization to the runnning time at the preceeding level of
discretization (where it is defined); and Predicted gives the
expected ratio of the running times (assuming $O(N^{1.5})$ asymptotic cost).

\label{section:deformed_torus}

\end{subsection}

%
%

\begin{subsection}{A heart-shaped domain.}

In these experiments, we solved the problem (\ref{formulations:neumann})
on the domain $\Omega$ shown in Figure~\ref{figure:heart}.
The boundary $\Gamma$ is described by the parameterization
\begin{equation*}
r(s,t) = \left(
\begin{array}{c}
\cos(s)\cos(t) + \cos^2(s) \sin^2(t) \\
\cos(s)\sin(t) \\
\sin(s)
\end{array}
\right),
\ \ \ 0\leq s,t< 2\pi,
\end{equation*}
although for the purposes of producing discretizations of the
surface $\Gamma$ a parameterization obtained
via projection onto a cube was used rather than this function $r(s,t)$.
The domain $\Omega$ is contained in the box with lower left corner
$(-1,-1,-1)$ and upper right corner $(1.25,1,1)$.

\begin{figure}[h!!]
\includegraphics[width=.3\textwidth]{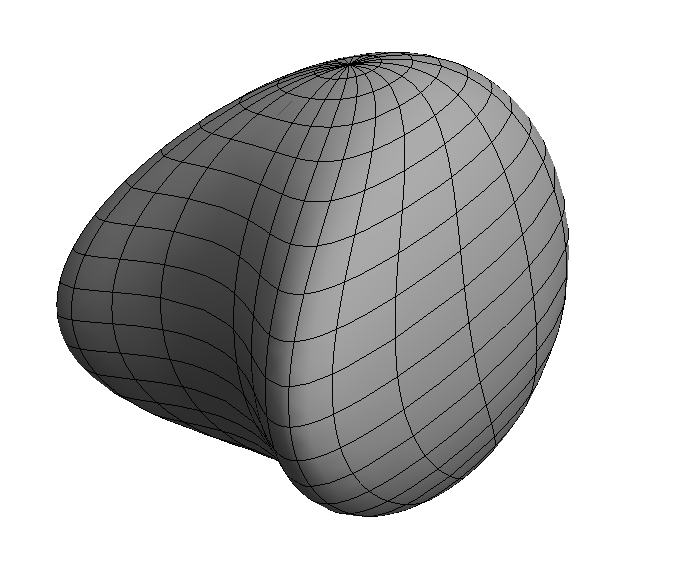}
\hfil
\includegraphics[width=.3\textwidth]{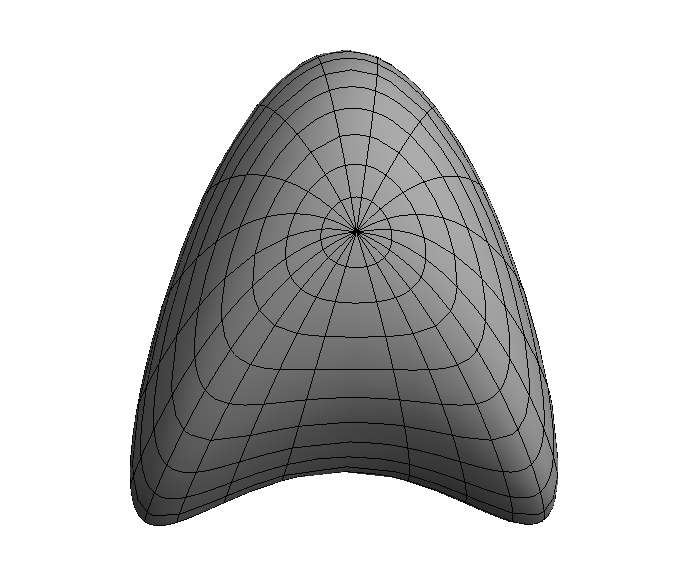}
\hfil
\includegraphics[width=.3\textwidth]{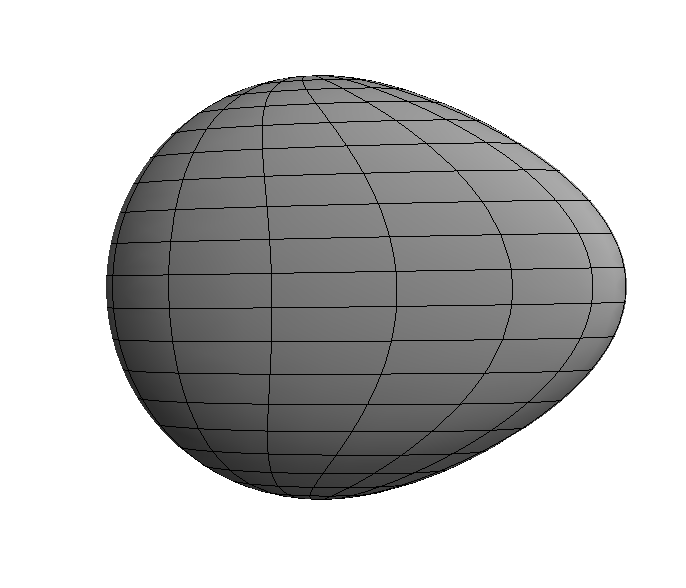}
\caption{Three views of the heart shaped surface of Section~\ref{section:heart_experiments}.}
\label{figure:heart}
\end{figure}

These experiments consisted of solving  the problem (\ref{formulations:neumann})
on the domain $\Omega$ at a variety of wavenumbers.  In each experiment,
a scattering matrix for the domain was formed and an instance
of the boundary value problem was solved using this scattering
matrix.  The boundary data was taken to be the derivative
of the restriction of the function
\begin{equation*}
g(x) = \frac{1}{4\pi} \frac{\exp(i \kappa |x|)}{|x|}
\end{equation*}
to $\Gamma$ with respect to the outward pointing unit normal
vector on $\Gamma$.  The error in the obtained charge distribution
was measured by computing
\begin{equation}
E =
 \left(\dblint_S \left| g(x)\right|^2 dx\right)^{-1/2}
\left(\dblint_S \left| g(x) - S_\kappa \sigma(x)\right|^2 dx\right)^{1/2},
\label{heart_experiments:error}
\end{equation}
where $S$ is the ball of radius $10$ centered at the point
$(-30,-30,-30)$ and $\sigma$ is the approximate charge distribution
obtained using the scattering matrix.
Discretizations of $\Gamma$ were obtained by triangulating
the faces of the cube which constitutes the parameterization domain.
An $8$th order discretization quadrature with $45$ points was applied
to each triangle.  The number
of triangles was increased with wavenumber.
In all cases, the desired precision for the scattering
matrix was taken to be  $\epsilon = 1.0\e{-3}$.
Table~\ref{table:heart} reports the number of discretization nodes $N$ on
the surface; the relative error $E$; the number of triangles $N_{\rm tris}$ into which
the parameterization domain was divided; the total running time
$T$; the dimensions $N_{\mbox{\tiny out}} \times N_{\mbox{\tiny in}}$
of the scattering matrix; the wavenumber $\kappa$ for the problem; and
the approximate diameter of the domain in wavelengths.

\begin{table}[h!!]
\small
\begin{tabular}{rrrrcccc}
\addlinespace[.5em]
$\kappa$ & Diameter & $N_{\mbox{\tiny tris}}$ & $N$ & $E$ & $T$ & $N_{\mbox{\tiny out}} \times N_{\mbox{\tiny in}}$ \\
\midrule
\addlinespace[.5em]
$2\pi$  & $3.46\lambda$  &     $48$ &   $2\,496$ & $4.20\e{-04}$ & $6.84\e{+00}$ & $251  \times 229$ \\
$4\pi$  & $6.92\lambda$  &    $192$ &   $9\,984$ & $1.33\e{-04}$ & $6.93\e{+01}$ & $571  \times 546$ \\
$6\pi$  & $10.39\lambda$ &    $768$ &  $39\,936$ & $1.88\e{-04}$ & $3.74\e{+02}$ & $996  \times 972$ \\
$8\pi$  & $13.87\lambda$ &    $768$ &  $39\,936$ & $3.61\e{-04}$ & $4.99\e{+02}$ & $1491 \times 1\,471$ \\
$10\pi$ & $17.32\lambda$ &   $3\,072$ & $159\,744$ & $4.96\e{-04}$ & $2.83\e{+03}$ & $2157 \times 2\,126$ \\
$12\pi$ & $20.79\lambda$ &   $3\,072$ & $159\,744$ & $1.89\e{-04}$ & $3.78\e{+03}$ & $2851 \times 2\,829$ \\
$14\pi$ & $24.25\lambda$ &  $12\,288$ & $638\,976$ & $6.75\e{-04}$ & $1.81\e{+04}$ & $3615 \times 3\,589$ \\
\end{tabular}
\vskip 1em
\caption{The results of the experiments of Section~\ref{section:heart_experiments}.}
\label{table:heart}
\vskip -2em
\end{table}

\label{section:heart_experiments}
\end{subsection}

\begin{subsection}{Scattering from a deformed cube.}
In this next set of experiments, we considered the problem
(\ref{formulations:neumann}) on a deformed cube $\Omega$
with $12$ curved edges and $8$ corner points.  The domain $\Omega$
is pictured in Figure~\ref{figure:cube}.  The top face of the boundary $\Gamma$
was parameterized over the region $-1 \leq s,t \leq 1$ via
\begin{equation*}
r(s,t) =
\left(2 + s^2 + t^2 \right)
\left(
\begin{array}{c}
s \\
t \\
1
\end{array}
\right)
\end{equation*}
and the other faces were handled in a like manner.  The domain $\Omega$ is
contained in the box with lower left corner $(-4,-4,-4)$
and upper right corner $(4,4,4)$.

\begin{figure}[h!!]
\includegraphics[width=.45\textwidth]{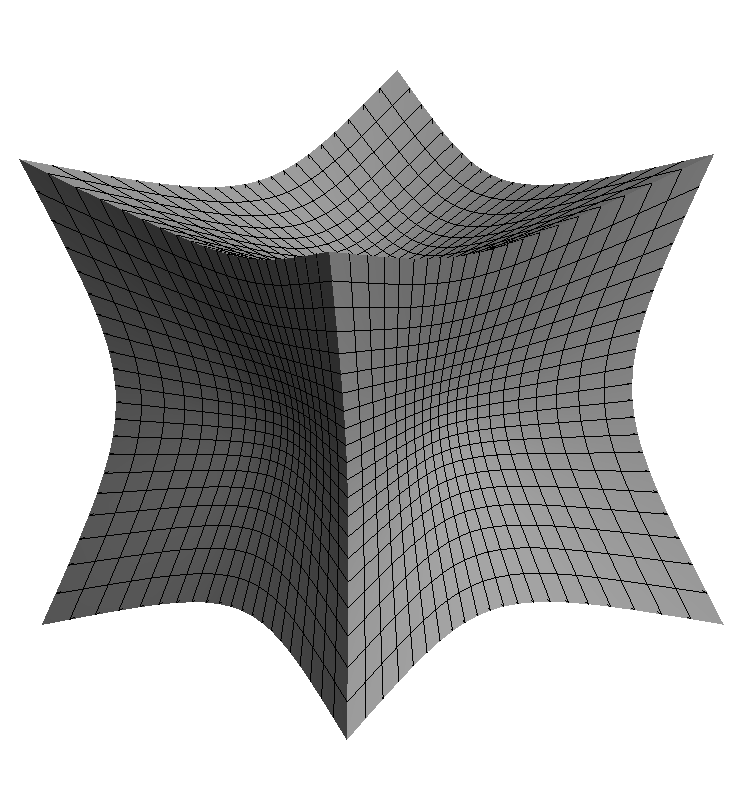}
\hfil
\includegraphics[width=.45\textwidth]{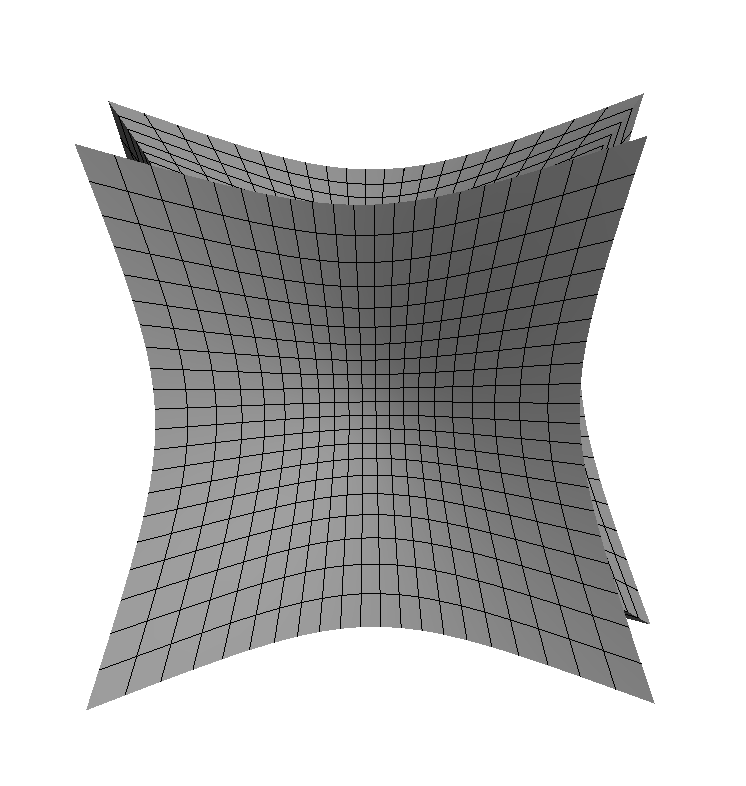}
\caption{The deformed cube domain considered in Section~\ref{section:cube_experiments}.
The domain has $12$ curved edges and $8$ corner points.
}
\label{figure:cube}
\end{figure}

\begin{figure}[h!!]
\includegraphics[width=.55\textwidth]{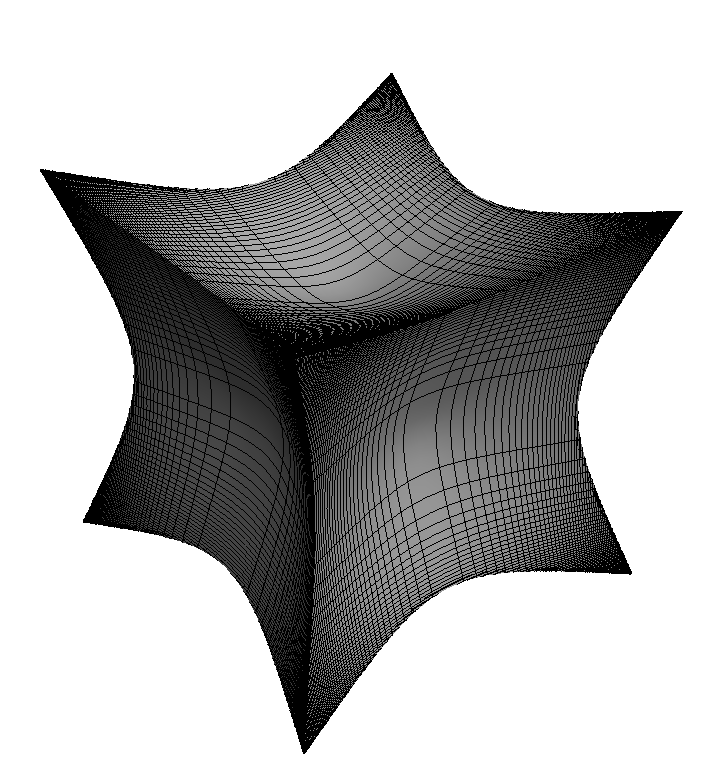}
\caption{Illustration of the unevenly refined grid used to generate discretization nodes on the deformed
cube surface considered in Section~\ref{section:cube_experiments}.}
\label{figure:cube2}
\end{figure}

In these experiments, the integral equation (\ref{formulations:naive}) defined
on $\Gamma$ was discretized by applying a $12$th order
discretization quadrature with $112$ nodes to increasing dense triangulations
of the parameterization domain.  We arranged for a denser distribution of discretization
nodes near the edges and corners of the deformed cube.
Figure~\ref{figure:cube2} shows one of the unevenly refined meshes used in these
experiments.

A scattering matrix was constructed for each discretization
and an instance of the boundary value problem was solved using this scattering
matrix.  In each experiment, the wavenumber was taken to be $\kappa = \pi/2$, making
 the domain roughly $3.46$ wavelengths in diameter, the desired precision
for the scattering matrix was set to $\epsilon=1.0\e{-10}$ and
the boundary data was taken to be the outward normal derivative
on $\Gamma$ of the function
\begin{equation*}
g(x) = \frac{1}{4\pi} \dblint_{\Gamma_{\mbox{\tiny in}}}
\frac{\exp(i \kappa |x-y|)}{|x-y|} ds(y),
\end{equation*}
where $\Gamma_{\mbox{\tiny in}}$ denotes the sphere of radius $0.1$ centered
at the point $(0.1,0,0)$. The error in the obtained charge distribution
$\sigma$ was measured via the quantity
\begin{equation}
E =
 \left(\dblint_{\Gamma_{\mbox{\tiny out}}} \left| g(x)\right|^2 dx\right)^{-1/2}
\left(\dblint_{\Gamma_{\mbox{\tiny out}}} \left| g(x) - S_\kappa \sigma(x)\right|^2 dx\right)^{1/2},
\label{cube_experiments:error}
\end{equation}
where $\sigma$ is the obtained charge distribution and
${\Gamma_{\mbox{\tiny out}}}$ denotes the sphere of radius $1$ centered
at the point $(-10,-10,0)$.
Table~\ref{table:cube} reports the error (\ref{cube_experiments:error}) along with $N_{\rm tris}$, $N$, $T$ and
$N_{\mbox{\tiny out}} \times N_{\mbox{\tiny in}}$ (whose meanings are as in the preceding sections).

\begin{table}[h!!]
\small
\begin{tabular}{rr@{\hspace{3em}}ccc}
\addlinespace[.5em]
 $N_{\mbox{\tiny tris}}$ & $N$ & $E$ & $T$ & $N_{\mbox{\tiny out}} \times N_{\mbox{\tiny in}}$\\
\midrule
\addlinespace[.5em]
$192$ & $21\,504$ & $2.60\e{-08}$ & $6.11\e{+02}$ &  $617 \times 712$\\
$432$ & $48\,384$ & $2.13\e{-09}$ & $1.65\e{+03}$ &  $620 \times 694$\\
$768$ & $86\,016$ & $3.13\e{-10}$ & $3.58\e{+03}$ &  $612 \times 685$\\
\end{tabular}
\vskip 1em
\caption{The results of the experiments preformed on the deformed cube of Section~\ref{section:cube_experiments}.}
\label{table:cube}
\vskip -2em
\end{table}

\label{section:cube_experiments}
\end{subsection}

\begin{subsection}{A Convergence Study involving a Trefoil Domain}
The experiments described in the preceding sections all involved
problems with artifical boundary data associated with known solutions.
In this section, we describe an experiment for which we did not have
access to an exact solution, and instead estimate the error by comparing
against a very finely resolved computed solution.

We considered the boundary value problem (\ref{formulations:dirichlet})
given on the trefoil domain $\Omega$ shown in Figure \ref{fig:trefoil}.
The parameterization of the boundary of $\Omega$ is given by
\begin{equation*}
r(s,t) = \vct{\gamma}(s) + \cos(s) \vct{n}(s) + \sin(-t) \vct{b}(s), \ \ \ 0\leq s,t < 2\pi,
\end{equation*}
where
\begin{equation*}
\vct{\gamma}(s)  =
\left(
\begin{array}{c}
\sin(3s)\\
\sin(s)+2\sin(2s)\\
\cos(s)-2\cos(s)
\end{array}
\right),
\end{equation*}
$\vct{n}(s)$ is the unit normal vector to the curve $\gamma(s)$ at $s$ and $\vct{b}(s)$ is the unit binormal
vector  at $s$.
The domain is contained in a box with bottom left corner at $[-1.5,-3, -3.5]$ and
top right corner $[1.5,3,3.5]$.

Discretizations of the surface
$\Gamma$ were formed by partitioning the parameter space into
$N_{\rm tris}$ triangles and utilizing an $8$th order (52 point) quadrature on each triangle
so  that $N = 52 N_{\rm tris}$.  A scattering matrix for each discretization
was then formed; $\Gamma_{\rm in}$
consisted of the point $(1,0,0)$ and $\Gamma_{\rm out}$ was taken to be
the surface of a sphere of radius $1$ centered at $(30,30,30)$.  The boundary
data was generated by a point charge on $\Gamma_{\rm in}$.  As the solution to
this problem is unknown, we report the $L^2$ convergence of the solution on
the surface of $\Gamma_{\rm out}$.  Specifically, we report
\begin{equation}
E_N= \left(\dblint_{\Gamma_{\rm out}} \left|(S_\kappa + i\kappa D_\kappa) \sigma_{4N}(x) - (S_\kappa + i\kappa D_\kappa) \sigma_N(x)\right|^2 dx\right)^{1/2},
\label{eq:conv_err}
\end{equation}
where $\sigma_{N}$ is the computed boundary charge distribution with $N$ discretization points.

We considered two choices of wavenumber and tolerance: (i) the domain was approximately $1\times 3\times 3$ wavelengths
in size and the tolerance is set to $\epsilon = 10^{-9}$; (ii) the domain was approximately
$2\times 6 \times 6$ wavelengths in size and the tolerance is set to $\epsilon = 10^{-6}$.
Tables \ref{tab:trefoilA} and \ref{tab:trefoilB} report on the convergence of the solution as the
number of discretization points $N$ is increased.

\begin{figure}[h!!]
\includegraphics[width=.45\textwidth]{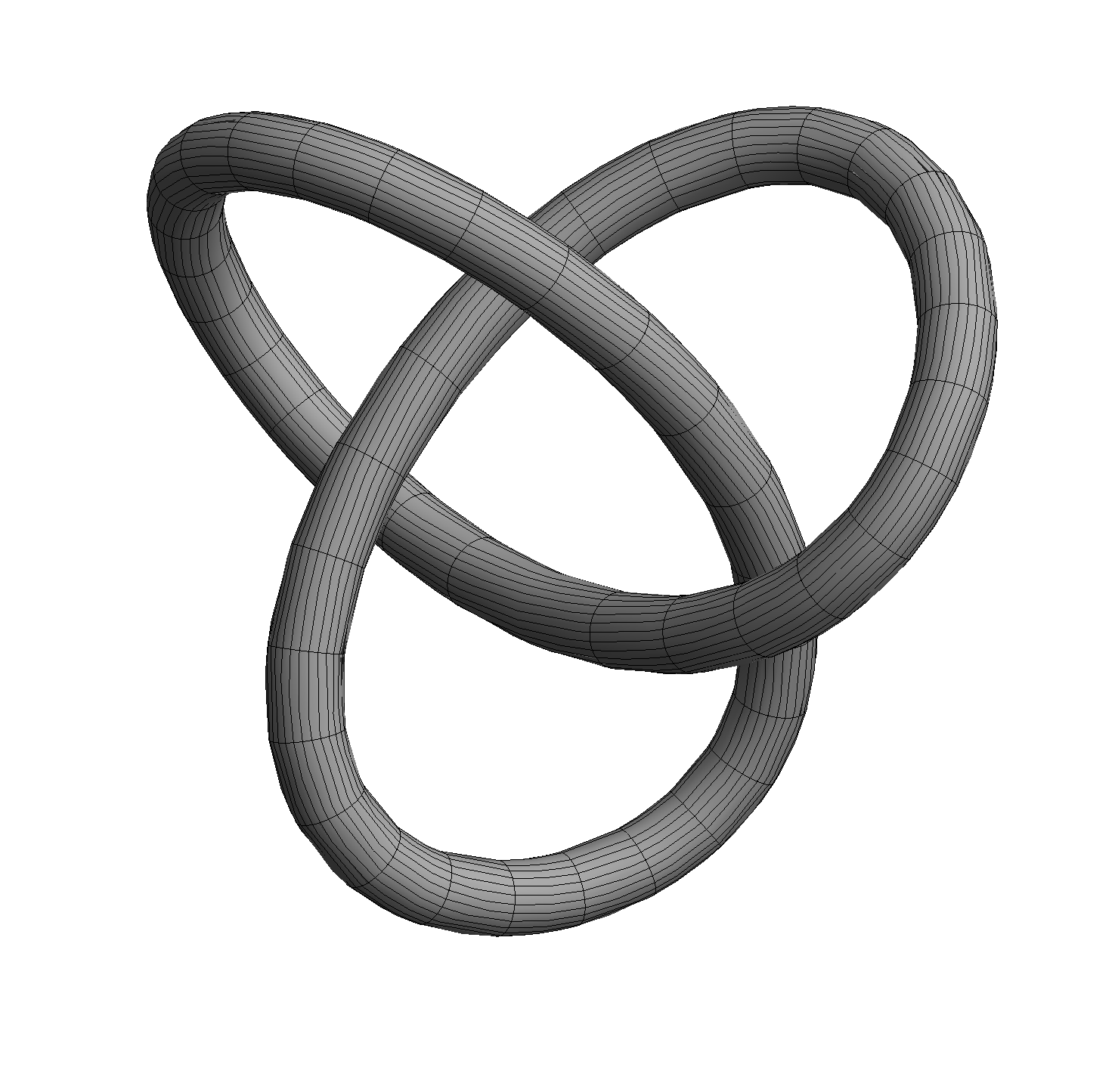}
\hfil
\includegraphics[width=.45\textwidth]{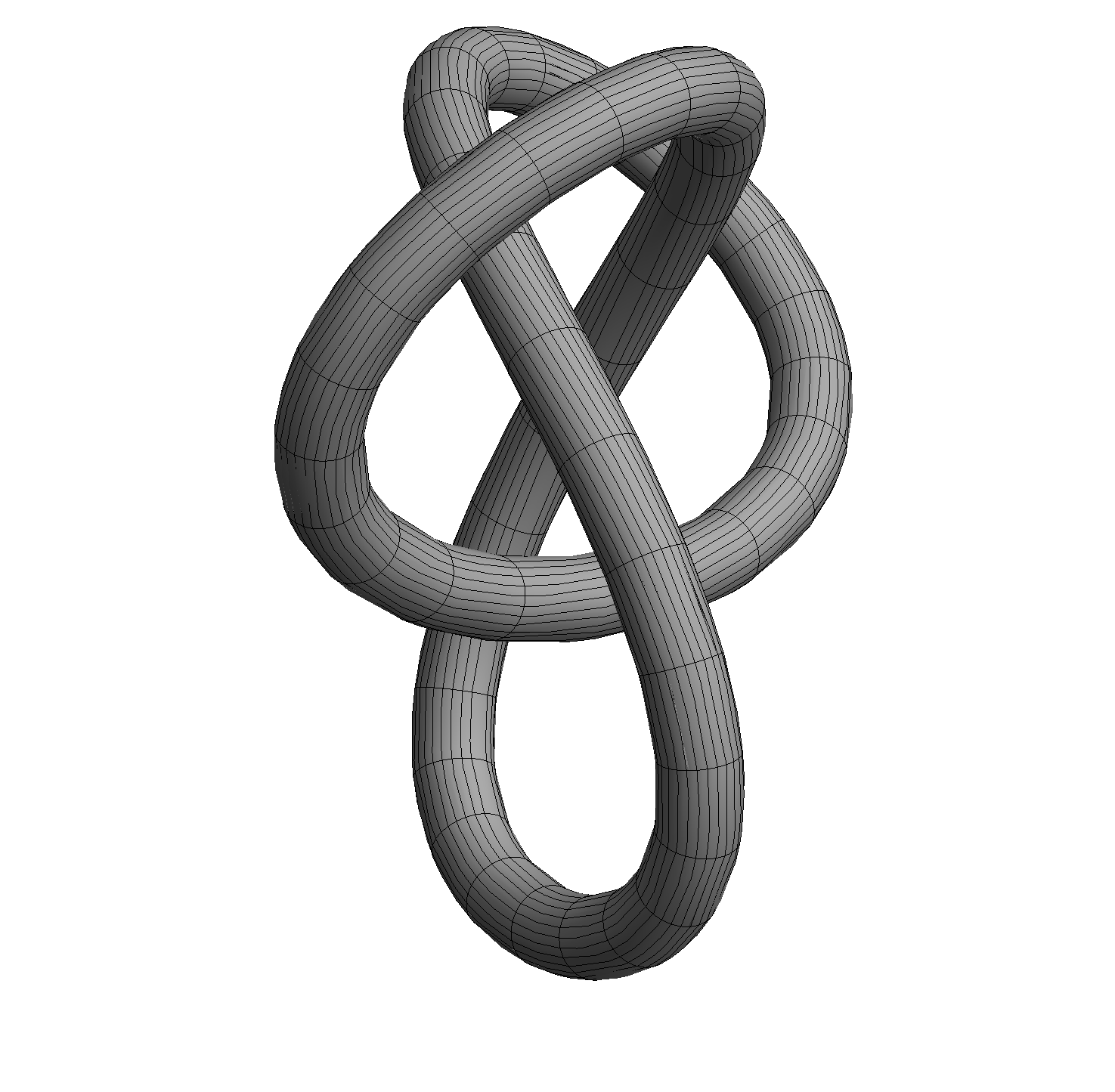}
\caption{The trefoil domain $\Omega$ considered in Section~\ref{section:trefoil_experiments}.
}
\label{fig:trefoil}
\end{figure}

\begin{table}[h!!]
\small
\begin{tabular}{rr@{\hspace{3em}}ccc}
\addlinespace[.5em]
 $N_{\mbox{\tiny tris}}$ & $N$ & $E_N$ & $T$ & $N_{\mbox{\tiny out}} \times N_{\mbox{\tiny in}}$\\
\midrule
\addlinespace[.5em]
$16$ & $832$ & $6.73\e{-04}$ & $1.17\e{+00}$ &  $754 \times 737$\\
$64$ & $3\,328$ & $2.33\e{-06}$ & $3.78\e{+01}$ &  $939 \times 910$\\
$256$ & $13\,312$ & $2.59\e{-08}$ & $3.61\e{+02}$ &  $945 \times 917$\\
$1\,024$ & $53\,248$ & $2.47\e{-11}$ & $2.55\e{+03}$ &  $948 \times 918$\\
$4\,096$ & $212\,992$ & - 		& $2.83\e{+04}$ &  $949 \times 921$\\
\end{tabular}
\vskip 1em
\caption{The results for the experiments performed on the trefoil of Section~\ref{section:trefoil_experiments} when
the domain is is approximately $1\times 3\times 3$ wavelengths
in size and $\epsilon = 10^{-9}$. }
\label{tab:trefoilA}
\vskip -2em
\end{table}

\begin{table}[h!!]
\small
\begin{tabular}{rr@{\hspace{3em}}ccc}
\addlinespace[.5em]
 $N_{\mbox{\tiny tris}}$ & $N$ & $E_N$ & $T$ & $N_{\mbox{\tiny out}} \times N_{\mbox{\tiny in}}$\\
\midrule
\addlinespace[.5em]
$16$&  $832 $ &      $2.87\e{-02}$ &     $1.94\e{+00}$ &  $832 \times 832$ \\
$64$&  $3\,328 $ &     $ 1.08\e{-03}$&    $7.59\e{+01}$&   $1\,771 \times 1\,745$\\
$256$&  $13\,312$ &    $5.57\e{-06 }$&    $ 4.12\e{+02}$& $1\,820 \times 1\,782$ \\
$1\,024$&  $53\,248 $ &   $ 3.04\e{-08}$& $2.16\e{+03}$&   $1\,822 \times 1\,788$\\
$4\,096$&  $212\,992$ &     -          & $1.27\e{+04}$ &  $1\,828 \times 1\,794$ \\
\end{tabular}
\vskip 1em
\caption{The results for the experiments performed on the trefoil of Section~\ref{section:trefoil_experiments} when
the domain is is approximately $2\times 6\times 6$ wavelengths
in size and $\epsilon = 10^{-6}$. }
\label{tab:trefoilB}
\vskip -2em
\end{table}

\label{section:trefoil_experiments}

\end{subsection}

\label{section:experiments}

\end{section}

\begin{section}{Acknowledgements}
J.~Bremer was supported by a fellowship from the Alfred
P. Sloan Foundation and by the Office of Naval Research under
contract N00014-12-1-0117.
P.G.~Martinsson was supported by the National Science Foundation under contracts
DMS0748488 and DMS0941476.

\end{section}

\begin{appendix}
\section{A robust integral formulation for solving (\ref{formulations:neumann})	 }
\label{app:nbie}

Consider the task of solving the Neumann boundary value problem (\ref{formulations:neumann}) using a
boundary integral equation formulation. In order to avoid the problem of spurious
resonances, it would be natural to use a combined layer representation of the form
\begin{equation}
\label{eq:tmpansatz}
u(x) = D_\kappa  \sigma(x) + i|\kappa|S_\kappa  \sigma(x).
\end{equation}
However, the BIE obtained by simply combining (\ref{eq:tmpansatz}) with the Neumann boundary
condition involves the \textit{hypersingular} operator
\begin{equation*}
N_\kappa \sigma(x) = \int_{\partial\Omega} \left( \eta_x \cdot \nabla_x \left( \eta_y \cdot \nabla_y G_\kappa (x,y) \right) \right) \sigma(y) ds(y).
\end{equation*}
Two difficulties are encountered when handling integral equations involving hypersingular operators:
\begin{enumerate}[1. ]

\item
Viewed as an operator from $L^2(\Gamma)$ to $L^2(\Gamma)$, $N_{\kappa}$ is not bounded.  In consequence,
if standard discretization methods are applied to integral equations involving this operator,
ill-conditioned linear systems will result.
This complication can be overcome by viewing $N_{\kappa}$ as on operator between two different
Sobolev spaces; however, more complicated analytic machinery would then be required.
\vskip 1em

\item
The singularity in the kernel of the hypersingular operator $N_\kappa$ is more severe than the singularities
of the standard operators of acoustic scattering.  As discussed in \cite{Bremer-Gimbutas2}, the approach
used here to evaluate singular integrals can be generalized to this case.  However, considerable
work is required and the resulting quadratures are somewhat less efficient than those used here.
\vskip 1em
\end{enumerate}

A technique known as ``regularization'' offers an alternative to handling the operator $N_\kappa$ \cite{Bruno-Elling-Turc}.
For  example, the representation (\ref{eq:combine}) can be replaced with
\begin{equation*}
u(x) = S_\kappa  \sigma(x) + D_\kappa  S_\kappa  \sigma(x)
\end{equation*}
which leads to the integral equation
\begin{equation}
 -\frac{1}{2}\sigma(x) + [N_\kappa  S_\kappa ] \sigma(x) +D_\kappa ^* \sigma(x) = g(x).
\label{eq:eq1}
\end{equation}
Using Calder\'on identities (which are discussed in, for instance, \cite{Nedelec}),  we can rewrite
(\ref{eq:eq1}) as
\begin{equation}
-\frac{3}{4} \sigma(x) + D_\kappa ^* \sigma(x) + (D_\kappa ^*)^2 \sigma(x) = g(x),
\label{formulations:composition}
\end{equation}
which does not involve hypersingular operators.
While the equation (\ref{formulations:composition}) is suitable for use
by iterative solvers, the composition of the operator $D_\kappa ^*$ makes
it cumbersome for direct solvers.

To make the formulation (\ref{formulations:composition}) more suitable for
direct solvers, we introduce an auxiliary variable
\begin{equation*}
\phi(x) = D_\kappa ^* \sigma(x),
\end{equation*}
which results to the system of integral equations
\begin{equation}
\begin{split}
-\frac{3}{4} \sigma(x) + D_\kappa ^* \sigma(x) + D_\kappa ^* \phi(x) &= g(x) \\
-\phi(x) + D_\kappa ^* \sigma(x) &= 0.
\end{split}
\label{formulation:system}
\end{equation}
The system (\ref{formulation:system}) is now ideal for
direct solvers as it does not involve hypersingular
operators and no compositions of operators appear.

\section{Derivation of scattering matrices and the merge formula from a field formulation}
\label{app:fields}

The derivation of the formula for merging two scattering matrices given in
Section \ref{sec:mergeformula} is quite formal. In this section, we provide
a different derivation that is hopefully more intuitive and gives a better idea
of the physical meaning of a scattering matrix $\mS_{\tau}$ for a subpatch
$\Gamma_{\tau} \subseteq \Gamma$. This derivation was omitted from the main
text since it requires the introduction of several auxiliary fields.

Throughout this section of the appendix, we rely on the hierarchical tree partitioning of
the domain $\Gamma$ described in Section \ref{sec:tree}, and the associated
HBS representation of a matrix given in Section \ref{sec:HBS}.

\subsection{Problem formulation} Our objective is to evaluate the map (\ref{eq:3step_map_disc}), which
is the discretization of the scattering problem (\ref{eq:3step_map}). We restate it for easy reference:
Given a vector $\vct{s}_{\rm in}$ representing a source distribution on the ``radiation source,'' we
seek to compute the field $\uu_{\rm out}$ on the ``antenna'' given by
\begin{equation}
\label{eq:3step_map_repeat}
\begin{array}{ccccccccc}
\uu_{\rm out} &=&
\mtx{A}_{\rm out} & \mtx{A}^{-1} & \mtx{A}_{\rm in} & \vct{s}_{\rm in}, \\[2mm]
N_{\rm out} \times 1 && N_{\rm out} \times N & N\times N & N\times N_{\rm in} & N_{\rm in} \times 1
\end{array}
\end{equation}
We define the \textit{external incoming field} on the scattering surface $\Gamma$ via
$$
\vct{u} = \mtx{A}_{\rm in }\vct{s}_{\rm in}.
$$
The core task is then to solve the equation
\begin{equation}
\label{eq:equil}
\mtx{A}\vct{\sigma} = \vct{u}
\end{equation}
to compute the field of induced charges $\vct{\sigma}$ on $\Gamma$. As described in Section
\ref{sec:disc_scatter_matrix}, we do not need to solve (\ref{eq:equil}) in full generality,
but only under two simplifying constraints:
\begin{enumerate}
\item The external incoming field $\vct{u}$ is restricted to a low-dimensional space spanned
by the provided basis matrix $\mhU_{1}$. In other words, there exists a short vector $\tuu_{1}$
such that $\vct{u} = \mhU_{1}\,\tuu_{1}$.
\item We only need the projection of $\vct{\sigma}$ onto a low-dimensional space spanned by
the provided basis matrix $\mhV_{1}$. In other words, we only need the short vector $\tilde{\vct{\sigma}}_{1} =
\mhV_{1}^{*}\vct{\sigma}$.
\end{enumerate}

For convenience in the hierarchical formulation of the scattering problem, we define for a
subpatch $\Gamma_{\tau}$ of the scattering surface an \textit{internal incoming field} via
\begin{equation}
\label{eq:kr1}
\vv_{\tau} = \mtx{A}(I_{\tau},I_{\tau}^{\rm c})\,\vct{\sigma}(I_{\tau}^{\rm c}).
\end{equation}
Then $\vv_{\tau}$ is the field on $\Gamma_{\tau}$ induced by charges on
$\Gamma \backslash \Gamma_{\tau}$. For a node $\tau$, the restriction of
(\ref{eq:equil}) to $I_{\tau}$,
\begin{equation}
\label{eq:kr2}
\mtx{A}(I_{\tau},I_{\tau})\,\vct{\sigma}(I_{\tau})
+
\mtx{A}(I_{\tau},I_{\tau}^{\rm c})\,\vct{\sigma}(I_{\tau}^{\rm c})
=
\uu(I_{\tau})
\end{equation}
can then be written
\begin{equation}
\label{eq:kr3}
\mtx{A}(I_{\tau},I_{\tau})\,\vct{\sigma}(I_{\tau})
=
\uu(I_{\tau}) - \vv_{\tau}.
\end{equation}
In words, equation (\ref{eq:kr3}) states that the induced charges $\vct{\sigma}(I_{\tau})$ on $\Gamma_{\tau}$
must create a local field that precisely match the \textit{total incoming field} $\uu(I_{\tau}) - \vv_{\tau}$,
which is the difference between the externally imposed incoming field $\uu(I_{\tau})$ and the internal
incoming field $\vv_{\tau}$ which is caused by induced charges on the rest of $\Gamma$.

\subsection{Short representations and the definition of the scattering matrix}
By definition, the columns of the long basis matrix $\mhU_{\tau}$ span both the
restriction $\uu(I_{\tau})$ of the external incoming field to the patch $\Gamma_{\tau}$,
and the columns of the off-diagonal block $\mtx{A}(I_{\tau},I_{\tau}^{\rm c})$.
As a consequence, there exist short vectors $\tuu_{\tau}$ and $\tvv_{\tau}$ such that
\begin{equation}
\label{eq:kr4}
\vv_{\tau} = \mhU_{\tau}\,\tvv_{\tau},
\qquad\mbox{and}\qquad
\uu(I_{\tau}) = \mhU_{\tau}\,\tuu_{\tau}.
\end{equation}
We say that $\tvv_{\tau}$ and $\tuu_{\tau}$ are the \textit{compressed representations}
of $\vv_{\tau}$ and $\uu(I_{\tau})$, respectively. Using these representations, (\ref{eq:kr3})
can be written
\begin{equation}
\label{eq:kr5}
\mtx{A}(I_{\tau},I_{\tau})\,\vct{\sigma}(I_{\tau})
=
\mhU_{\tau}\,\bigl(\tuu_{\tau} - \tvv_{\tau}\bigr).
\end{equation}
We also define a compressed representation of the induced charges $\vct{\sigma}(I_{\tau})$ on $\Gamma_{\tau}$ via
\begin{equation}
\label{eq:kr6}
\tilde{\vct{\sigma}}_{\tau}
=
\mhV_{\tau}^{*}
\vct{\sigma}(I_{\tau}).
\end{equation}
Multiplying (\ref{eq:kr5}) by $\mhV_{\tau}^{*}\mA(I_{\tau},I_{\tau})^{-1}$, we now get the scattering equation
\begin{equation}
\label{eq:meaningS}
\tilde{\vct{\sigma}}_{\tau}
=
\mS_{\tau}\,\bigl(\tuu_{\tau} - \tvv_{\tau}\bigr),
\end{equation}
where, recall from (\ref{eq:def_S1}),
\begin{equation}
\label{eq:kr8}
\mS_{\tau} = \mhV_{\tau}^{*}\mA(I_{\tau},I_{\tau})^{-1}\mhU_{\tau}.
\end{equation}
Equation (\ref{eq:meaningS}) states that \textit{the scattering
matrix $\mS_{\tau}$ maps the short representation
$\tuu_{\tau} - \tvv_{\tau}$ of the total incoming field on $\Gamma_{\tau}$
to the short representation $\tilde{\vct{\sigma}}_{\tau}$ of the induced charges.}

\subsection{Hierarchical construction of scattering matrices}
For leaf nodes $\tau$, we compute $\mS_{\tau}$ directly from the definition (\ref{eq:kr8}).
For parent nodes, we compute the scattering matrices via a hierarchical merging procedure;
it turns out that if we know the scattering matrices $\mtx{S}_{\alpha}$ and $\mtx{S}_{\beta}$
of the children $\{\alpha,\beta\}$ of a node $\tau$, then $\mtx{S}_{\tau}$ can be built
inexpensively as follows: Observe that (\ref{eq:kr5}) can be written
\begin{equation}
\label{eq:kr9}
\left[\begin{array}{cc}
\mtx{A}(I_{\alpha},I_{\alpha}) & \mtx{A}(I_{\alpha},I_{\beta})\\
\mtx{A}(I_{\beta},I_{\alpha})  & \mtx{A}(I_{\beta},I_{\beta})
\end{array}\right]
\left[\begin{array}{c}\vct{\sigma}(I_{\alpha}) \\ \vct{\sigma}(I_{\beta}) \end{array}\right]
=
\mhU_{\tau}\,\bigl(\tuu_{\tau} - \tvv_{\tau}\bigr).
\end{equation}
Inserting into (\ref{eq:kr9}) the factorizations (\ref{eq:sibling_factorization})
of $\mtx{A}(I_{\alpha},I_{\beta})$ and $\mtx{A}(I_{\beta},I_{\alpha})$, and the hierarchical representation
(\ref{eq:hbases_repeat}) of $\mhU_{\tau}$, we get
\begin{equation}
\label{eq:kr10}
\left[\begin{array}{cc}
\mtx{A}(I_{\alpha},I_{\alpha}) & \mhU_{\alpha}\mtA_{\alpha,\beta}\mhV_{\beta}^{*}\\
\mhU_{\beta}\mtA_{\beta,\alpha}\mhV_{\alpha}^{*} & \mtx{A}(I_{\beta},I_{\beta})
\end{array}\right]
\left[\begin{array}{c}\vct{\sigma}(I_{\alpha}) \\ \vct{\sigma}(I_{\beta}) \end{array}\right]
=
\left[\begin{array}{cc}
\mhU_{\alpha} & \mtx{0} \\
\mtx{0} & \mhU_{\beta}
\end{array}\right]\,
\mU_{\tau}\,\bigl(\tuu_{\tau} - \tvv_{\tau}\bigr).
\end{equation}
Left-multiply (\ref{eq:kr10}) by $\mtwo{\mhV_{\alpha}^{*}\mtx{A}(I_{\alpha},I_{\alpha})^{-1}}{\mtx{0}}{\mtx{0}}{\mhV_{\beta}^{*}\mtx{A}(I_{\beta},I_{\beta})^{-1}}$
to obtain
\begin{equation}
\label{eq:kr11}
\left[\begin{array}{cc}
\mtx{I} & \mS_{\alpha}\mtA_{\alpha,\beta} \\
\mS_{\beta}\mtA_{\beta,\alpha} & \mtx{I}
\end{array}\right]
\left[\begin{array}{c} \tilde{\vct{\sigma}}_{\alpha} \\ \tilde{\vct{\sigma}}_{\beta}\end{array}\right]
=
\left[\begin{array}{cc}
\mS_{\alpha} & \mtx{0} \\ \mtx{0} & \mS_{\beta}
\end{array}\right]\,
\mU_{\tau}\,\bigl(\tuu_{\tau} - \tvv_{\tau}\bigr).
\end{equation}
Solving (\ref{eq:kr11}) for $\{\tilde{\vct{\sigma}}_{\alpha},\tilde{\vct{\sigma}}_{\beta}\}$ and left multiplying the
result by $\mV_{\tau}^{*}$, we find
\begin{equation}
\label{eq:kr12}
\underbrace{\mV_{\tau}^{*}
\left[\begin{array}{c} \tilde{\vct{\sigma}}_{\alpha} \\ \tilde{\vct{\sigma}}_{\beta}\end{array}\right]}_{=\tilde{\vct{\sigma}}_{\tau}}
=
\underbrace{\mV_{\tau}^{*}
\left[\begin{array}{cc}
\mtx{I} & \mS_{\alpha}\mtA_{\alpha,\beta} \\
\mS_{\beta}\mtA_{\beta,\alpha} & \mtx{I}
\end{array}\right]^{-1}
\left[\begin{array}{cc}
\mS_{\alpha} & \mtx{0} \\ \mtx{0} & \mS_{\beta}
\end{array}\right]\,
\mU_{\tau}}_{=\mS_{\tau}}
\bigl(\tuu_{\tau} - \tvv_{\tau}\bigr).
\end{equation}
Comparing equation (\ref{eq:kr12}) to (\ref{eq:meaningS}), we rediscover
the formula (\ref{eq:mergelemmaA}) for the merge of two scattering matrices
given in Lemma \ref{lemma:mergelemma}.

\end{appendix}

\bibliographystyle{acm}
\bibliography{surface,main_bib}

\end{document}